\documentclass[a4paper,twoside]{amsart}

\usepackage[pdftex]{color,graphicx} 
\usepackage{amssymb,amsfonts,amsmath}

\usepackage{multirow}		

\usepackage{tikz}
\usetikzlibrary{matrix,arrows}
\usetikzlibrary{positioning}
\usetikzlibrary{fit}
\usetikzlibrary{patterns}

\usepackage{stmaryrd}		

\def\lqq{``\hspace{0.15pt}\nolinebreak}
\def\rqq{\hspace{0.5pt}\nolinebreak''}

\newcommand{\cF}{{\mathcal F}}
\newcommand{\cG}{{\mathcal G}}

\newcommand{\tom}{\varnothing}					

\newcommand{\CC}{\mathbb{C}}							

\newcommand{\sep}{\,|\,} 								

\newcommand{\dbrac}[1]{{\llbracket #1 \rrbracket}} 	

\newcommand{\Fl}{\mathrm{Flags}}						

\newcommand{\rma}{\mathrm{a}}							
\newcommand{\rmc}{\mathrm{c}}	
\newcommand{\rmi}{\mathrm{i}}	
\newcommand{\rmr}{\mathrm{r}}	

\newcommand{\brur}[2]{(#1\! \leftrightarrow\! #2)}		

\newcommand{\vinc}[3]{								
\begin{tikzpicture}[baseline = (X.base)]
	\useasboundingbox (0.1,0) rectangle (#1*0.23,0.1);
	\foreach \x/\y in {#2}
	{
		\draw (\x*0.2,0) node (X) {$\y$};
	}
	
	\foreach \z in {#3}
	{
		\ifnum 0<\z
			\ifnum \z<#1
				\draw[thick] (\z*0.2-0.07,-0.17) -- (\z*0.2+0.27,-0.17);
			\fi
		\fi
		
		\ifnum 0=\z
			\draw[thick] (0.1,0.1) -- (0.1,-0.17) -- (0.21,-0.17);
		\fi
		
		\ifnum \z=#1
			\draw[thick] (\z*0.2+0.11,0.1) -- (\z*0.2+0.11,-0.17) -- (\z*0.2,-0.17);
		\fi
	}
\end{tikzpicture}
}

\newcommand{\vinci}[3]{\!\!							
\begin{tikzpicture}[baseline = (X.base)]
	\useasboundingbox (0.1,0) rectangle (#1*0.17,0.16);

	\foreach \x/\y in {#2}
	{
		\draw (\x*0.15,0) node (X) {$\scriptstyle\y$};
	}
	
	\foreach \z in {#3}
	{
		\ifnum 0<\z
			\ifnum \z<#1
				\draw[thick] (\z*0.15-0.06,-0.06) -- (\z*0.15+0.2,-0.06);
			\fi
		\fi
		
		\ifnum 0=\z
			\draw[thick] (0.1,0.1) -- (0.1,-0.17) -- (0.21,-0.17);
		\fi
		
		\ifnum \z=#1
			\draw[thick] (\z*0.2+0.11,0.1) -- (\z*0.2+0.11,-0.17) -- (\z*0.2,-0.17);
		\fi
	}
	\useasboundingbox (0,0);
\end{tikzpicture}
	\!
}

\newcommand{\bivinc}[4]{\!							
\begin{tikzpicture}[baseline]

	\foreach \x/\y in {#2}
	{
		\draw (\x*0.2,0.3) node {$\x$};
		\draw (\x*0.2,0) node {$\y$};
	}
	
	\foreach \z in {#3}
	{
		\ifnum 0<\z
			\ifnum \z<#1
				\draw[thick] (\z*0.2-0.07,-0.17) -- (\z*0.2+0.27,-0.17);
			\fi
		\fi
		
		\ifnum 0=\z
			\draw[thick] (0.1,0.1) -- (0.1,-0.17) -- (0.21,-0.17);
		\fi
		
		\ifnum \z=#1
			\draw[thick] (\z*0.2+0.11,0.1) -- (\z*0.2+0.11,-0.17) -- (\z*0.2,-0.17);
		\fi
	}
	
	\foreach \z in {#4}
	{
		\ifnum 0<\z
			\ifnum \z<#1
				\draw[thick] (\z*0.2-0.07,0.47) -- (\z*0.2+0.27,0.47);
			\fi
		\fi
		
		\ifnum 0=\z
			\draw[thick] (0.1,0.21) -- (0.1,0.47) -- (0.21,0.47);
		\fi
		
		\ifnum \z=#1
			\draw[thick] (\z*0.2+0.11,0.21) -- (\z*0.2+0.11,0.47) -- (\z*0.2,0.47);
		\fi
	}
\end{tikzpicture}
\!
}

\newcommand{\bivincs}[4]{\!										
\begin{tikzpicture}[baseline]

	\foreach \x/\w/\y in {#2}
	{
		\draw (\x*0.2,0.3) node {$\w$};
		\draw (\x*0.2,0) node {$\y$};
	}
	
	\foreach \z in {#3}
	{
		\ifnum 0<\z
			\ifnum \z<#1
				\draw[thick] (\z*0.2-0.07,-0.17) -- (\z*0.2+0.27,-0.17);
			\fi
		\fi
		
		\ifnum 0=\z
			\draw[thick] (0.1,0.1) -- (0.1,-0.17) -- (0.21,-0.17);
		\fi
		
		\ifnum \z=#1
			\draw[thick] (\z*0.2+0.11,0.1) -- (\z*0.2+0.11,-0.17) -- (\z*0.2,-0.17);
		\fi
	}
	
	\foreach \z in {#4}
	{
		\ifnum 0<\z
			\ifnum \z<#1
				\draw[thick] (\z*0.2-0.07,0.47) -- (\z*0.2+0.27,0.47);
			\fi
		\fi
		
		\ifnum 0=\z
			\draw[thick] (0.1,0.21) -- (0.1,0.47) -- (0.21,0.47);
		\fi
		
		\ifnum \z=#1
			\draw[thick] (\z*0.2+0.11,0.21) -- (\z*0.2+0.11,0.47) -- (\z*0.2,0.47);
		\fi
	}
\end{tikzpicture}
\!
}

\newcommand{\bivinci}[4]{\!										
\begin{tikzpicture}[baseline]
	\useasboundingbox (0.1,0) rectangle (#1*0.17,0.5);

	\foreach \x/\y in {#2}
	{
		\draw (\x*0.15,0.225) node {$\scriptstyle\x$};
		\draw (\x*0.15,0) node {$\scriptstyle\y$};
	}
	
	\foreach \z in {#3}
	{
		\ifnum 0<\z
			\ifnum \z<#1
				\draw[thick] (\z*0.15-0.06,-0.06) -- (\z*0.15+0.2,-0.06);
			\fi
		\fi
		
		\ifnum 0=\z
			\draw[thick] (0.1,0.1) -- (0.1,-0.17) -- (0.21,-0.17);
		\fi
		
		\ifnum \z=#1
			\draw[thick] (\z*0.2+0.11,0.1) -- (\z*0.2+0.11,-0.17) -- (\z*0.2,-0.17);
		\fi
	}
	
	\foreach \z in {#4}
	{
		\ifnum 0<\z
			\ifnum \z<#1
				\draw[thick] (\z*0.15-0.06,0.43) -- (\z*0.15+0.2,0.43);
			\fi
		\fi
		
		\ifnum 0=\z
			\draw[thick] (0.1,0.21) -- (0.1,0.47) -- (0.21,0.47);
		\fi
		
		\ifnum \z=#1
			\draw[thick] (\z*0.2+0.11,0.21) -- (\z*0.2+0.11,0.47) -- (\z*0.2,0.47);
		\fi
	}
\end{tikzpicture}
\!
}

\newcommand{\pattern}[4]{										
  \raisebox{0.6ex}{
  \begin{tikzpicture}[scale=0.35, baseline=(current bounding box.center), #1]
    \foreach \x/\y in {#4}
      \fill[pattern=north east lines] (\x,\y) rectangle +(1,1);
    \draw (0.01,0.01) grid (#2+0.99,#2+0.99);
    \foreach \x/\y in {#3}
      \filldraw (\x,\y) circle (6pt);
  \end{tikzpicture}}
}

\newcommand{\impattern}[5]{									
  \raisebox{0.6ex}{
  \begin{tikzpicture}[scale=0.35, baseline=(current bounding box.center), #1]
    \foreach \x/\y in {#5}
      \fill[pattern=north east lines] (\x,\y) rectangle +(1,1);
    \draw (0.01,0.01) grid (#2+0.99,#2+0.99);
    \foreach \x/\y in {#4}
      \draw[fill=white] (\x,\y) circle (6pt);
    \foreach \x/\y in {#3}
      \filldraw (\x,\y) circle (6pt);
  \end{tikzpicture}}
}

\newcommand{\imopattern}[6]{									
  \raisebox{0.6ex}{
  \begin{tikzpicture}[scale=0.35, baseline=(current bounding box.center), #1]
    \foreach \x/\y in {#6}
      \fill[pattern=north east lines] (\x,\y) rectangle +(1,1);
    \draw (0.01,0.01) grid (#2+0.99,#2+0.99);
    \foreach \x/\y in {#4}
      \draw[fill=white] (\x,\y) circle (6pt);
    \foreach \x/\y in {#5}
      \draw[fill=white] (\x,\y) circle (10pt);
    \foreach \x/\y in {#3}
      \filldraw (\x,\y) circle (6pt);
  \end{tikzpicture}}
}

\newcommand{\patternsbm}[5]{									
  \raisebox{0.6ex}{
  \begin{tikzpicture}[scale=0.35, baseline=(current bounding box.center), #1]
    \foreach \x/\y in {#4}
      \fill[pattern=north east lines] (\x,\y) rectangle +(1,1);
    \draw (0.01,0.01) grid (#2+0.99,#2+0.99);
    \foreach \x/\y/\z/\w/\A in {#5}
       \fill[color = gray!100, opacity=0.5, rounded corners] (\x+0.075,\y+0.075) rectangle (\z-0.075,\w-0.075);
    \foreach \x/\y/\z/\w/\A in {#5}
       \draw[color = black, rounded corners] (\x+0.075,\y+0.075) rectangle (\z-0.075,\w-0.075);
    \foreach \x/\y/\z/\w/\A in {#5}
       \fill[black] (\x/2+\z/2,\y/2+\w/2) node {$\scriptstyle\A$};
    \foreach \x/\y in {#3}
      \filldraw (\x,\y) circle (6pt);

  \end{tikzpicture}}
}

\newcommand{\patternsbmm}[6]{									
  \raisebox{0.6ex}{											
  \begin{tikzpicture}[scale=0.35, baseline=(current bounding box.center), #1]
    \foreach \x/\y in {#4}
      \fill[pattern=north east lines] (\x,\y) rectangle +(1,1);
    \draw (0.01,0.01) grid (#2+0.99,#2+0.99);
    \foreach \x/\y/\z/\w/\A in {#5}
       \fill[color = gray!100, opacity=0.5, rounded corners] (\x+0.075,\y+0.075) rectangle (\z-0.075,\w-0.075);
    \foreach \x/\y/\z/\w/\A in {#5}
       \draw[color = black, rounded corners] (\x+0.075,\y+0.075) rectangle (\z-0.075,\w-0.075);
    \foreach \x/\y/\z/\w/\A in {#6}
       \fill[black] (\x/2+\z/2,\y/2+\w/2) node {$\scriptstyle\A$};
    \foreach \x/\y in {#3}
      \filldraw (\x,\y) circle (6pt);

  \end{tikzpicture}}
}

\newcommand{\patternsbmfreelybraided}[6]{						
  \raisebox{0.6ex}{
  \begin{tikzpicture}[scale=0.35, baseline=(current bounding box.center), #1]
    \foreach \x/\y in {#4}
      \fill[pattern=north east lines] (\x,\y) rectangle +(1,1);
    \draw (0.01,0.01) grid (#2+0.99,#2+0.99);
    \foreach \x/\y/\z/\w/\A in {#5}
       \fill[color = gray!100, opacity=0.5, rounded corners] (\x+0.075,\y+0.075) rectangle (\z-0.075,\w-0.075);
    \foreach \x/\y/\z/\w/\A in {#5}
       \draw[color = black, rounded corners] (\x+0.075,\y+0.075) rectangle (\z-0.075,\w-0.075);
    \fill[color = gray!100, opacity=0.5, rounded corners] (0+0.075,3) -- (0+0.075,2+0.075) -- (1-0.075,2+0.075) -- (1-0.1,3+0.1) --
    (2-0.075,3+0.075) -- (2-0.075,4-0.075) -- (0+0.075,4-0.075) -- (0+0.075,3);
    \draw[color = black, rounded corners] (0+0.075,3) -- (0+0.075,2+0.075) -- (1-0.075,2+0.075) -- (1-0.1,3+0.1) --
    (2-0.075,3+0.075) -- (2-0.075,4-0.075) -- (0+0.075,4-0.075) -- (0+0.075,3);
    \fill[color = gray!100, opacity=0.5, rounded corners] (4-0.075,1) -- (4-0.075,2-0.075) -- (3+0.075,2-0.075) -- (3+0.1,1-0.1) --
    (2+0.075,1-0.075) -- (2+0.075,0+0.075) -- (4-0.075,0+0.075) -- (4-0.075,1);
    \draw[color = black, rounded corners] (4-0.075,1) -- (4-0.075,2-0.075) -- (3+0.075,2-0.075) -- (3+0.1,1-0.1) --
    (2+0.075,1-0.075) -- (2+0.075,0+0.075) -- (4-0.075,0+0.075) -- (4-0.075,1);
    \foreach \x/\y/\z/\w/\A in {#6}
       \fill[black] (\x/2+\z/2,\y/2+\w/2) node {$\scriptstyle\A$};
    \foreach \x/\y in {#3}
      \filldraw (\x,\y) circle (6pt);

  \end{tikzpicture}}
}

{  \theoremstyle{plain} 
   \newtheorem*{named-thm}{}
   \newtheorem*{fullyr�ing}{Fullyr�ing}
   \newtheorem{theorem}{Theorem}
   \newtheorem{proposition}[theorem]{Proposition} 
   \newtheorem{lemma}[theorem]{Lemma}
   \newtheorem{corollary}[theorem]{Corollary}
   
   \newtheorem{example}[theorem]{Example}

}
{\theoremstyle{definition}
   \newtheorem{definition}[theorem]{Definition}
   \newtheorem{defn-prop}[theorem]{Skilgreining-Fullyr�ing}

}

\usepackage[pdftex]{hyperref}

\begin{document}

\title{A unification of permutation patterns related to Schubert varieties
}


\author{Henning \'Ulfarsson
}


\address{School of Computer Science,
Reykjavik University, Reykjav\'ik, Iceland}


\begin{abstract}
 We obtain new connections between permutation patterns and singularities
 of Schubert varieties, by giving a new characterization of Gorenstein
 varieties in terms of so called \emph{bivincular patterns}.
 These are generalizations of classical patterns where conditions are placed on the
 location of an occurrence in a permutation, as well as on the values in
 the occurrence.
 This clarifies what happens when the requirement of smoothness is weakened
 to factoriality and further to Gorensteinness, extending work of
 Bousquet-M\'elou and Butler (2007), and Woo and Yong (2006).
 We also show how \emph{mesh patterns}, introduced by Br\"and\'en and Claesson (2011),
 subsume many other types of patterns and define an extension of them
 called \emph{marked mesh patterns}. We use these new patterns to further simplify
 the description of Gorenstein Schubert varieties and give a new description of Schubert
 varieties that are \emph{defined by inclusions}, introduced by Gasharov and Reiner (2002).
 We also give a description of $123$-hexagon avoiding permutations, introduced by Billey and Warrington (2001),
 Dumont permutations and cycles in terms of marked mesh patterns.
\end{abstract}

\maketitle


\section{Introduction}

In this paper we exhibit new connections between permutation patterns and singularities
of Schubert varieties $X_\pi$ in the complete flag variety $\Fl(\CC^n)$,
by giving a new characterization of Gorenstein varieties in terms of which \emph{bivincular patterns} the
permutation $\pi$ avoids. Bivincular patterns, defined by Bousquet-M\'elou, Claesson,
Dukes and Kitaev~\cite{BCDK10},
are generalizations of classical patterns where conditions are placed on the
location of an occurrence in a permutation, as well as on the values in
the occurrence.
This clarifies what happens when the requirement of smoothness is weakened
to factoriality and further to Gorensteinness, extending work of
Bousquet-M\'elou and Butler~\cite{MR2376109}, and Woo and Yong~\cite{MR2264071}.
We also prove results that translate some known patterns in the literature into
bivincular patterns. In particular we will give a characterization of the
Baxter permutations.

Table \ref{table} summarizes the main results in the paper related to bivincular patterns.
The first line in the table is due to Ryan~\cite{MR870962}, Wolper~\cite{MR1013667} and Lakshmibai and Sandhya~\cite{MR1051089}
and says that a Schubert variety $X_\pi$ is non-singular (or smooth) if and only if $\pi$ avoids
the patterns $1324$ and $2143$. Note that some authors use a different convention for the correspondence between permutations
and Schubert varieties, which results in the reversal of the permutations. These authors would then use the patterns $4231$, $3412$ to
identify the smooth Schubert varieties. Saying that the variety
$X_\pi$ is non-singular means that every local ring is regular.

\begin{table}[ht]
\caption{Connections between singularity properties and bivincular patterns
}
\label{table}
\begin{tabular}{ l l l }
\hline\noalign{\smallskip}
$X_\pi$ is	 	& \multicolumn{2}{l}{The permutation $\pi$ avoids the patterns} \\
\noalign{\smallskip}\hline\noalign{\smallskip}
smooth			& $\vinc{4}{1/2, 2/1, 3/4, 4/3}{}$ and $\vinc{4}{1/1, 2/3, 3/2, 4/4}{}$ & \\
\noalign{\smallskip}\noalign{\smallskip}\noalign{\smallskip}
factorial		& $\vinc{4}{1/2, 2/1, 3/4, 4/3}{2}$ and $\vinc{4}{1/1, 2/3, 3/2, 4/4}{}$ & \\
\noalign{\smallskip}\noalign{\smallskip}\noalign{\smallskip}
\multirow{2}{*}{Gorenstein}		& \multirow{2}{*}{$\bivinc{5}{1/3,2/1,3/5,4/2,5/4}{2}{3}$,
				  $\bivinc{5}{1/2,2/4,3/1,4/5,5/3}{3}{2}$} &
				  and associated Grassmannians avoid \\
				  & & two bivincular pattern families \\
\noalign{\smallskip}\hline
\end{tabular}
\end{table}%

A weakening of
this condition is the requirement that every local ring only be a unique
factorization domain; a variety satisfying this is a \emph{factorial} variety.
Bousquet-M\'elou and Butler~\cite{MR2376109} proved a conjecture
stated by Woo and Yong (personal communication) that factorial
Schubert varieties are those that correspond to permutations avoiding
$1324$ and bar-avoiding $21\overline{3}54$. In
the terminology of Woo and Yong~\cite{MR2264071} the bar-avoidance of the latter
pattern corresponds to avoiding $2143$ with Bruhat restriction
$\brur{1}{4}$, or equivalently, interval avoiding
$[2413, 2143]$ in the terminology of Woo and Yong~\cite{MR2422304}.  However, as remarked by
Steingr\'imsson~\cite{steingrimsson-2008}, bar-avoiding $21\overline{3}54$
is equivalent to avoiding the vincular pattern
$\vinc{4}{1/2, 2/1, 3/4, 4/3}{2}$. See Theorem \ref{thm:Ulfarsson-Factorial} for the details.

A further weakening is to only require that the local rings of $X_\pi$ be
Gorenstein local rings, in which case we say that $X_\pi$ is a \emph{Gorenstein}
variety. Woo and Yong~\cite{MR2264071} showed that $X_\pi$ is Gorenstein
if and only if it avoids two patterns with two Bruhat restrictions each,
as well as satisfying a certain condition on descents. We will translate
their results into avoidance of bivincular patterns;
see Theorem \ref{thm:Ulfarsson-Gorenstein}.

We also show how \emph{mesh patterns}, introduced by Br\"and\'en and Claesson~\cite{BC11},
subsume many other types of patterns, such as interval patterns defined by Woo and
Yong~\cite{MR2422304}, and define an extension of them
called \emph{marked mesh patterns}. We use these new patterns to further simplify
the description of Gorenstein Schubert varieties (see Theorem \ref{thm:Ulfarsson-Gorenstein2})
and give a new description of Schubert
varieties that are \emph{defined by inclusions}, introduced by Gasharov and Reiner~\cite{MR1934291}
(see Theorem \ref{thm:dbi}).
We also give a description of $123$-hexagon avoiding permutations, introduced by Billey and Warrington~\cite{MR1826948},
in terms of the avoidance of $123$ and one marked mesh pattern (see Proposition \ref{prop:123hex}).
Finally, in Example \ref{ex:markedmesh}, we describe Dumont permutations~\cite{0297.05004} (of the first and second kind) and cycles with marked mesh patterns.

\begin{table}[ht]
\caption{Connections between singularity properties and marked mesh patterns}
\label{table2}
\begin{tabular}{ l l l }
\hline\noalign{\smallskip}
$X_\pi$ is	 	& \multicolumn{2}{l}{The permutation $\pi$ avoids the patterns} \\
\noalign{\smallskip}\hline\noalign{\smallskip}
smooth			& \multicolumn{2}{l}{$\patternsbm{scale=0.75}{ 4 }{ 1/2, 2/1, 3/4, 4/3 }{}{}$ and $\vinc{4}{1/1, 2/3, 3/2, 4/4}{}$} \\
\noalign{\smallskip}\noalign{\smallskip}\noalign{\smallskip}
factorial		& \multicolumn{2}{l}{$\patternsbm{scale=0.75}{ 4 }{ 1/2, 2/1, 3/4, 4/3 }{2/0, 2/1, 2/2, 2/3, 2/4}{}$ and $\vinc{4}{1/1, 2/3, 3/2, 4/4}{}$} \\
\noalign{\smallskip}\noalign{\smallskip}\noalign{\smallskip}
defined by inclusions &
\multicolumn{2}{l}{$\patternsbm{scale=0.75}{ 4 }{ 1/2, 2/1, 3/4, 4/3 }{}{3/1/4/2/{},1/3/2/4/\scriptscriptstyle{1}}$,
$\patternsbm{scale=0.75}{ 4 }{ 1/2, 2/1, 3/4, 4/3 }{}{4/1/5/2/\scriptscriptstyle{1},3/0/4/1/\scriptscriptstyle{1}}$ and $\vinc{4}{1/1, 2/3, 3/2, 4/4}{}$} \\
\noalign{\smallskip}\noalign{\smallskip}\noalign{\smallskip}
\multirow{2}{*}{Gorenstein}	& \multirow{2}{*}{$\patternsbm{scale=0.75}{ 4 }{ 1/2, 2/1, 3/4, 4/3 }{2/0, 2/1, 2/2, 2/3, 2/4, 0/2, 1/2, 3/2, 4/2}{3/1/4/2/{},1/3/2/4/\scriptscriptstyle{1}}$} &
				  and associated Grassmannians avoid \\
				  & & two mesh pattern families \\
\noalign{\smallskip}\noalign{\smallskip}\noalign{\smallskip}
\noalign{\smallskip}\noalign{\smallskip}\noalign{\smallskip}
\noalign{\smallskip}\noalign{\smallskip}\noalign{\smallskip}
$123$-hexagon av.\ & \!$\patternsbm{scale=0.75}{ 4 }{ 1/2, 2/1, 3/4, 4/3 }{}{2/4/3/5/\scriptscriptstyle{1}, 0/2/1/3/\scriptscriptstyle{1}, 4/2/5/3/\scriptscriptstyle{1}, 2/0/3/1/\scriptscriptstyle{1}}$ \\
\noalign{\smallskip}\hline
\end{tabular}
\end{table}%

\section{Three types of pattern avoidance} \label{sec:three-types}
Here we recall definitions of different types of patterns. We will
use one-line notation for all permutations, \textit{e.g.}, write $\pi = 312$ for
the permutation in $S_3$ that satisfies $\pi(1) = 3$, $\pi(2) = 1$ and
$\pi(3) = 2$.

The three types correspond to:
\begin{itemize}
 
 \item Bivincular patterns, subsuming vincular patterns and classical patterns.
 
 \item Barred patterns.
 
 \item Bruhat-restricted patterns.
 
\end{itemize}

\subsection{Bivincular patterns}
We denote the symmetric group on $n$ letters by $S_n$ and refer to its
elements as \emph{permutations}. We write permutations as words $\pi
= a_1a_2\dotsm a_n$, where the letters are distinct and come from the
set $\{1, 2, \dotsc, n\}$.  A \emph{pattern} $p$ is also a
permutation, but we are interested in when a pattern is
\emph{contained} in a permutation $\pi$ as described below.

An \emph{occurrence} (or \emph{embedding}) of a pattern $p$ in a
permutation $\pi$ is classically defined as a subsequence in $\pi$, of
the same length as $p$, whose letters are in the same relative order
(with respect to size) as those in $p$. For example, the pattern $123$
corresponds to an increasing subsequence of three letters in a
permutation. If we use the notation $1_\pi$ to denote the first,
$2_\pi$ for the second and $3_\pi$ for the third letter in an occurrence,
then we are simply requiring that $1_\pi < 2_\pi < 3_\pi$.
If a permutation has no occurrence of a pattern $p$ we
say that $\pi$ \emph{avoids} $p$.

\begin{example}
 The permutation $32415$ contains two occurrences of the pattern $123$
 corresponding to the sub-words $345$ and $245$. It avoids the pattern
 $132$.
\end{example}
The occurrence of a pattern in a permutation $\pi$ can also be defined
as a subset of the diagram $G(p) = \{ (i,\pi(i)) \sep 1 \leq i \leq n \}$,
that \lqq looks like\rqq\ the diagram of the pattern. Below is the diagram of
the pattern $123$ and two copies of the digram of the permutation $32415$
where we have indicated the two occurrences of the pattern by circling
the dots.
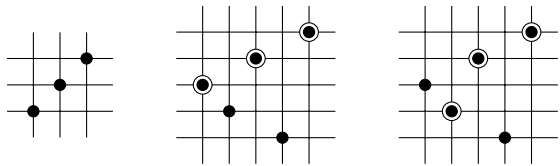
\[
\pattern{scale=1}{3}{1/1, 2/2, 3/3}{}
\qquad
\imopattern{scale=1}{5}{1/3, 2/2, 3/4, 4/1, 5/5}{}{1/3, 3/4, 5/5}{}
\qquad
\imopattern{scale=1}{5}{1/3, 2/2, 3/4, 4/1, 5/5}{}{2/2, 3/4, 5/5}{}
\]
In a \emph{vincular pattern} two adjacent letters may or may not
be underlined. If they are underlined it means that the
corresponding letters in the permutation $\pi$ must be adjacent.

\begin{example}
 The permutation $32415$ contains one occurrence of the pattern
 $\vinc{3}{1/1,2/2,3/3}{1}$ corresponding to the sub-word $245$. It avoids the 
 pattern $\vinc{3}{1/1,2/2,3/3}{2}$.
 The permutation $\pi = 324615$ has one occurrence of the
 pattern $2143$, namely the sub-word $3265$, but no
 occurrence of $\vinc{4}{1/2, 2/1, 3/4, 4/3}{2}$, since $2$ and $6$ are not adjacent in
 $\pi$.
\end{example}

It is also convenient to consider vincular patterns as certain types of diagrams.
We use dark vertical strips between dots that are required to be adjacent in the
pattern. Notice that only the second occurrence of the classical pattern $123$ satisfies
the requirements of the vincular pattern, since in the former the dot corresponding
to $2$ in $\pi$ lies in the forbidden strip.
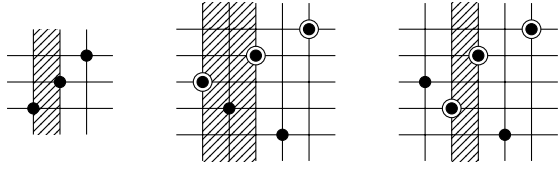
\[
\pattern{scale=1}{3}{1/1, 2/2, 3/3}{1/0, 1/1, 1/2, 1/3}
\qquad
\imopattern{scale=1}{5}{1/3, 2/2, 3/4, 4/1, 5/5}{}{1/3, 3/4, 5/5}{1/0, 1/1, 1/2, 1/3, 1/4, 1/5, 2/0, 2/1, 2/2, 2/3, 2/4, 2/5}
\qquad
\imopattern{scale=1}{5}{1/3, 2/2, 3/4, 4/1, 5/5}{}{2/2, 3/4, 5/5}{2/0, 2/1, 2/2, 2/3, 2/4, 2/5}
\]
\noindent
These types of patterns have been studied sporadically for
a very long time but were not defined in full generality until
Babson and Steingr\'imsson~\cite{MR1758852}.

This notion was generalized further in Bousquet-M\'elou et al~\cite{BCDK10}:
In a \emph{bivincular pattern} we are also allowed to place
restrictions on the values that occur in an embedding of a pattern. We
use two-line notation to describe these patterns. If there is a line over
the letters $i$, $i+1$ in the top row, it means that the corresponding letters in an occurrence
must be adjacent in values. This
is best described by an example: 

\begin{example}
 An occurrence of the pattern
 $\bivinc{3}{1/1,2/2,3/3}{}{2}$ in a permutation $\pi$ is an increasing subsequence
 of three letters, such that the third one is larger than the second by exactly
 $1$, or more simply, $3_\pi = 2_\pi + 1$.
 The permutation $32415$ contains two occurrence of this bivincular pattern
 corresponding to the sub-words $345$ and $245$. The second one is
 also an occurrence of $\bivinc{3}{1/1,2/2,3/3}{1}{2}$. The permutation
 avoids the bivincular pattern $\bivinc{3}{1/1,2/2,3/3}{2}{1,2}$.
\end{example}

By also using horizontal strips we are able to draw diagrams of bivincular patterns.
Below is the diagram of $\bivinc{3}{1/1,2/2,3/3}{1}{2}$ together with one occurrence
of it.
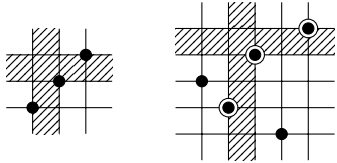
\[
\pattern{scale=1}{3}{1/1, 2/2, 3/3}{1/0, 1/1, 1/2, 1/3, 0/2, 2/2, 3/2}
\qquad
\imopattern{scale=1}{5}{1/3, 2/2, 3/4, 4/1, 5/5}{}{2/2, 3/4, 5/5}{2/0, 2/1, 2/2, 2/3, 2/4, 2/5,
0/4, 1/4, 3/4, 4/4, 5/4}
\]

\noindent
We will also use the notation of \cite{BCDK10} to write bivincular patterns:
A bivincular pattern consists of a triple $(p,X,Y)$ where $p$ is a permutation
in $S_k$ and $X,Y$ are subsets of $\dbrac{0,k} = \{0,\dotsc,k\}$. With this notation an
occurrence of a bivincular pattern in a permutation $\pi = \pi_1\dotsm\pi_n$ in $S_n$
is a subsequence $\pi_{i_1}\dotsm\pi_{i_k}$ such that the letters in the subsequence are in
the same relative order as the letters of $p$ and
\begin{itemize}

 \item for all $x$ in $X$, $i_{x+1} = i_x + 1$; and
 
 \item for all $y$ in $Y$, $j_{y+1} = j_y + 1$, where
 $\{ \pi_{i_1}, \dotsc, \pi_{i_k} \} = \{ j_1, \dotsc, j_k \}$ and
 $j_1 < j_2 < \dotsm < j_k$.
 
\end{itemize}
By convention we put $i_0 = 0 = j_0$ and $i_{k+1} = n+1 = j_{k+1}$.

\begin{example}
 We can translate all of the patterns we have discussed above into this
 notation:
  \begin{align*}
 	123 									&= (123,\tom,\tom),
	&132 								&= (132,\tom,\tom),
	&\vinc{3}{1/1,2/2,3/3}{1}				&= (123,\{1\},\tom),\\
	\vinc{3}{1/1,2/2,3/3}{2}				&= (123,\{2\},\tom),
	&2143				 				&= (2143,\tom,\tom),
	&\vinc{4}{1/2, 2/1, 3/4, 4/3}{2}		&= (2143,\{2\},\tom),\\
	\bivinc{3}{1/1,2/2,3/3}{}{1} 			&= (123,\tom,\{1\}),
	&\bivinc{3}{1/1,2/2,3/3}{}{1,2}	 	&= (123,\tom,\{1,2\}),
	&\bivinc{3}{1/1,2/2,3/3}{2}{1,2}		&= (123,\{2\},\{1,2\}).
 \end{align*}
\end{example}

\noindent
We have not considered the case when $0$ or $k$ are elements of $X$ or $Y$,
as we will not need those cases. We just remark that if $0 \in X$ then
an occurrence of $(p,X,Y)$ must start at the beginning of a permutation $\pi$,
in other words, $\pi_{i_1} = \pi_1$. The other cases are similar.

The bivincular patterns behave well with respect to the operations 
reverse, complement and inverse: Given a bivincular pattern $(p,X,Y)$
we define
\begin{align*}
	(p,X,Y)^\rmr &= (p^\rmr,k-X,Y), \qquad
	(p,X,Y)^\rmc = (p^\rmc,X,k-Y), \\
	(p,X,Y)^\rmi &= (p^\rmi,Y,X),
\end{align*}
where $p^\rmr$ is the usual reverse of the permutation $p$,
$p^\rmc$ is the usual complement of the permutation $p$, and
$p^\rmi$ is the usual inverse of the permutation $p$. Here
$k-M = \{ k-m \sep m \in M\}$. In \cite{BCDK10} the following is proved.
\begin{lemma}
 Let $\rma$ denote one of the operations above, or a composition of them.
 Then a permutation $\pi$ avoids the bivincular pattern $p$ if and only
 if the permutation $\pi^\rma$ avoids the bivincular pattern
 $p^\rma$.
\end{lemma}

\subsection{Barred patterns}
We will only consider a single pattern of this type, but the general definition
is easily inferred from this special case.
We say that a permutation $\pi$ \emph{avoids the barred pattern}
$21\overline{3}54$ if $\pi$ avoids the pattern
$2143$ (corresponding to the unbarred elements)
\emph{except} where that pattern is a part of the pattern
$21354$.  This notation for barred patterns was
introduced by West~\cite{W90}. It turns out that avoiding this barred
pattern is equivalent to avoiding the vincular pattern $\vinc{4}{1/2, 2/1, 3/4, 4/3}{2}$; see section
\ref{sec:connections}. See also section \ref{sec:intmesh} on how to write barred
patterns as mesh patterns.

\begin{example}
 The permutation $\pi = 4257613$ avoids the barred pattern
 $21\overline{3}54$ since
 the unique occurrence of $2143$, as the sub-word
 $4276$, is contained in the sub-word $42576$ which is an occurrence of
 $21354$. Note that it also avoids $\vinc{4}{1/2, 2/1, 3/4, 4/3}{2}$.
\end{example}

\subsection{Bruhat-restricted patterns}
We recall the definition of Bruhat-restricted patterns from
Woo and Yong~\cite{MR2264071}. First we need the \emph{Bruhat order}
on permutations in $S_n$, defined as follows: Given integers $i<j$ in
$\dbrac{1,n}$ and a permutation $\pi \in S_n$ we define $\pi\brur{i}{j}$
as the permutation that we get from $\pi$ by swapping $\pi(i)$ and $\pi(j)$.
For example $24153\brur{1}{4} = 54123$. We then say that $\pi\brur{i}{j}$
\emph{covers} $\pi$ if $\pi(i) < \pi(j)$ and for every $k$ with $i < k < j$
we have either $\pi(k) < \pi(i)$ or $\pi(k) > \pi(j)$. We then define the
Bruhat order as the transitive closure of the above covering relation.
This definition should be compared to the construction of the graph
$G_\pi$ in subsection \ref{subsec:FacSchubs-forest-like-perms}.
We see that in our example above that $24153\brur{1}{4}$ does not
cover $24153$ since we have $\pi(2) = 4$. Now, given a pattern
$p$ with a set of transpositions $\mathcal{T} = \{\brur{i_\ell}{j_\ell}\}$
we say that a permutation $\pi$ \emph{contains $(p,\mathcal{T})$}, or that
\emph{$\pi$ contains the Bruhat-restricted pattern $p$} (if $\mathcal{T}$
is understood from the context), if there is an embedding of $p$ in $\pi$
such that if any of the transpositions in $\mathcal{T}$ are carried out
on the embedding the resulting permutation covers $\pi$.

We should note that Bruhat-restricted patterns were further generalized
to \emph{intervals of patterns} in Woo and Yong~\cite{MR2422304}.
We delay the discussion of this type of pattern avoidance until
section \ref{sec:intmesh}, where we also introduce \emph{mesh patterns}
and show that an interval pattern is a special case of a mesh pattern.

In the next section we will show how these three types of patterns are
related to one another.

\section{Connections between the three types} \label{sec:connections}

\subsection{Factorial Schubert varieties and forest-like permutations}
\label{subsec:FacSchubs-forest-like-perms}

Bousquet-M\'elou and Butler~\cite{MR2376109} defined and studied
\emph{forest-like} permutations. Here we recall their
definition: Given a permutation
$\pi$ in $S_n$, construct a graph $G_\pi$ on the vertex set $\{1, 2, \dotsc, n\}$
by joining $i$ and $j$ if
\begin{enumerate}
 
 \item $i < j$ and $\pi(i) < \pi(j)$; and

 \item there is no $k$ such that $i < k < j$ and $\pi(i) < \pi(k) < \pi(j)$.
 
\end{enumerate}
\noindent
The permutation $\pi$ is \emph{forest-like} if the graph $G_\pi$ is a forest.
In light of the definition of Bruhat covering above we see that the vertices
$i$ and $j$ are connected in the graph of $G_\pi$ if and only if
$\pi\brur{i}{j}$ covers $\pi$.
 
They then show that a permutation is forest-like if and only if it
avoids the classical pattern $1324$ and the barred pattern
$p_\text{bar}�= 21\overline{3}54$. This barred pattern can
be described in terms of Bruhat-restricted embeddings and in terms of bivincular
patterns, as we now show.

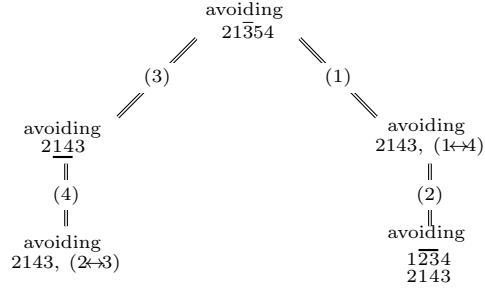
\begin{figure}
 \begin{tikzpicture}[description/.style={fill=white,inner sep=2pt}]
 
 \matrix (m) [matrix of math nodes, row sep=3em, column sep=2.5em,
 text height=1.5ex, text depth=0.25ex]
 { & \substack{\text{avoiding } \\ 21\overline{3}54} &   \\
 \substack{\text{avoiding } \\ \vinci{4}{1/2, 2/1, 3/4, 4/3}{2}} &
 & \substack{\text{avoiding } \\ 2143,\ \brur{1}{4}} \\
 \substack{\text{avoiding } \\ 2143,\ \brur{2}{3}} &
 & \substack{\text{avoiding } \\ \bivinci{4}{1/2,2/1,3/4,4/3}{}{2}}�\\
  };
 
 \path[-,font=\scriptsize]
 (m-1-2.south west) edge [double] node[description]
 {(\ref{enum:barred-Bruhat-bivincular3})}
 (m-2-1.north east)
 
 (m-1-2.south east) edge [double] node[description]
 {(\ref{enum:barred-Bruhat-bivincular1})}
 (m-2-3.160)
 
 (m-2-1) edge [double, shorten >=3pt, shorten <= 3pt] node[description]
 {(\ref{enum:barred-Bruhat-bivincular4})}
 (m-3-1)
 
 (m-2-3) edge [double, shorten >=3pt, shorten <= 3pt] node[description]
 {(\ref{enum:barred-Bruhat-bivincular2})}
 (m-3-3);
 
 \end{tikzpicture}
 \caption{The barred pattern $21\overline{3}54$ gives a
 connection between two bivincular patterns. The labels on the edges
 correspond to the enumerated list below}
 \label{fig:patterns}
\end{figure}
 
\begin{enumerate}
 
 \item \label{enum:barred-Bruhat-bivincular1}
 Bousquet-M\'elou and Butler \cite{MR2376109} remark that forest-like
 permutations $\pi$ correspond to factorial Schubert varieties $X_\pi$ and
 avoiding the barred pattern is equivalent to avoiding
 $p_\text{Br} = 2143$ with Bruhat restriction $\brur{1}{4}$.
 This last part is easily verified.
 
 \item \label{enum:barred-Bruhat-bivincular2}
 Avoiding $p_\text{Br} = 2143$ with Bruhat restriction $\brur{1}{4}$
 is equivalent to avoiding the bivincular pattern
 $p_\text{bi} = \bivinc{4}{1/2,2/1,3/4,4/3}{}{2}$,
 as we will now show:
 
 Assume $\pi$ contains the bivincular pattern $p_\text{bi}$,
 so we can find an embedding of it in $\pi$ such that $3_\pi = 2_\pi + 1$. This
 embedding clearly satisfies the Bruhat restriction.
 
 Now assume that $\pi$ has an embedding of $p_\text{Br}$. If $3_\pi = 2_\pi + 1$ we
 are done. Otherwise $2_\pi + 1$ is either to the right of $3_\pi$ or to the left
 of $2_\pi$ (because of the Bruhat restriction). In the first
 case change $3_\pi$ to $2_\pi + 1$ and we are done. In the second case replace
 $2_\pi$ with $2_\pi + 1$, thus reducing the distance in values to $3_\pi$,
 then repeat.
 
 \item \label{enum:barred-Bruhat-bivincular3}
 The barred pattern $p_\text{bar}�= 21\overline{3}54$ has
 another connection to bivincular patterns: avoiding it is equivalent
 to avoiding the bivincular pattern $q_\text{biv} = \vinc{4}{1/2, 2/1, 3/4, 4/3}{2}$,
 as remarked in the survey by Steingr\'imsson~\cite{steingrimsson-2008}.
 
 \item \label{enum:barred-Bruhat-bivincular4}
 We can translate this into Bruhat-restricted embeddings as well: Avoiding
 the bivincular pattern $q_\text{bi} = \vinc{4}{1/2, 2/1, 3/4, 4/3}{2}$ is equivalent to
 avoiding $q_\text{Br} = 2143$ with Bruhat restriction $\brur{2}{3}$:
 
 Assume $\pi$ has an embedding of $q_\text{Br}$. If $1_\pi$ and $4_\pi$ are
 adjacent then we are done. Otherwise look at the letter to right of $1_\pi$.
 If this letter is larger than $4_w$ we can replace $4_w$ by it and we are done.
 Otherwise this letter must be less than $4_w$, which implies by the Bruhat
 restriction, that it must also be less than $1_w$. In this case we replace
 $1_w$ by this letter, and repeat.
 
 Now assume $\pi$ has an embedding of the bivincular pattern $q_\text{bi}$.
 If $1_\pi$ and $4_\pi$ are adjacent we are done. Otherwise look at the letter to the
 right of $1_\pi$. This letter is either smaller than $1_\pi$ or larger than $4_\pi$.
 In the first case, replace $1_\pi$ with this letter; in the second case,
 replace $4_\pi$ with this letter. Then repeat if necessary.
 
\end{enumerate}

The above argument will be generalized in Proposition \ref{prop:Bruhat-bivinc},
but this special case gives us:

\begin{theorem}
\label{thm:Ulfarsson-Factorial}[\cite{MR2376109},\cite{steingrimsson-2008}]
 Let $\pi \in S_n$. The Schubert variety $X_\pi$ is factorial if and only if
 $\pi$ avoids the patterns $\vinc{4}{1/2, 2/1, 3/4, 4/3}{2}$ and $1324$.
\end{theorem}
 
From the equivalence of the patterns in Figure \ref{fig:patterns}
we also get that a permutation $\pi$ avoids the bivincular pattern
\[
	\vinc{4}{1/2, 2/1, 3/4, 4/3}{2} = (2143,\{2\},\tom)
\]
if and only if it avoids
\[
 \bivinc{4}{1/2, 2/1, 3/4, 4/3}{}{2} = (2143,\tom,\{2\}).
\]
We will prove this without going through the barred pattern, and
then generalize the proof, but first of all we should note that these
bivincular patterns are inverses of one another, and that will simplify
the proof.

 Assume $\pi$ contains $\bivinc{4}{1/2, 2/1, 3/4, 4/3}{}{2}$. If $1_\pi$ and $4_\pi$
 are adjacent in $\pi$ we are done. Otherwise consider the element immediately
 to the right of $1_\pi$. If this element is less than $2_\pi$ then replace $1_\pi$
 by it
 and we will have reduced the distance between $1_\pi$ and $4_\pi$. If this element
 is larger than $2_\pi$ it must also be larger than $3_\pi$, since
 $3_\pi = 2_\pi + 1$,
 so replace $4_\pi$ by it. This will (immediately, or after several steps) produce
 an occurrence of $\vinc{4}{1/2, 2/1, 3/4, 4/3}{2}$.
 
 Now assume $\pi$ contains $\vinc{4}{1/2, 2/1, 3/4, 4/3}{2}$. Then
 $\pi^\rmi$ contains the inverse pattern
 \[
 	(\vinc{4}{1/2, 2/1, 3/4, 4/3}{2})^\rmi = \bivinc{4}{1/2, 2/1, 3/4, 4/3}{}{2}.
 \]
 Then by the above, $\pi^\rmi$ contains
 $\vinc{4}{1/2, 2/1, 3/4, 4/3}{2}$, so $\pi = (\pi^\rmi)^\rmi$ contains
 $(\vinc{4}{1/2, 2/1, 3/4, 4/3}{2})^\rmi = 
 \bivinc{4}{1/2, 2/1, 3/4, 4/3}{}{2}$.
 
This generalizes to:
 
\begin{proposition} \label{prop:bivinc-bivinc}
 Let $p$ be the pattern
 \[
 	\vinc{8}{1/\cdot, 2/\cdot, 3/\cdot, 4/1, 5/k, 6/\cdot, 7/\cdot, 8/\cdot}{4}
	=
	(\dotsm 1k \dotsm, \{j\},\tom)
 \]
 in $S_k$, where $j = p^\rmi(1)$ is the \emph{index} of $1$ in $p$.
 A permutation $\pi$ in $S_n$ that avoids the pattern $p$ must also avoid the
 bivincular pattern
 \[	
 	\bivincs{8}{1/1/\cdot, 2/2/\cdot, 3/\cdot/\cdot, 4/\cdot/1,
	            5/\cdot/k, 6/\cdot/\cdot, 7/\cdot/\cdot, 8/k/\cdot}{}{2,3,4,5,6}=
	            (\dotsm 1 k \dotsm, \tom, \{2,3,\dotsc,k-2\}).
 \]
\end{proposition}

\begin{proof}
 Assume a permutation $\pi$ contains the latter pattern in the proposition.
 If $1_\pi$ and $k_\pi$ are adjacent in $\pi$ we are done. Otherwise
 consider the element immediately to the right of $1_\pi$. If this element
 is larger than $(k-1)_\pi$ we replace $k_\pi$ by it and are done. Otherwise
 this element must be less than $(k-1)_\pi$ and therefore less than $2_\pi$,
 so we can replace $1_\pi$ by it, and repeat.
\end{proof}

By applying the reverse to everything in Proposition
\ref{prop:bivinc-bivinc} we get:

\begin{corollary}
 Let $p$ be the pattern
 \[
 	\vinc{8}{1/\cdot, 2/\cdot, 3/\cdot, 4/k, 5/1, 6/\cdot, 7/\cdot, 8/\cdot}{4}
	=
	(\dotsm k 1 \dotsm, \{j\},\tom)
 \]
 in $S_k$, where $j = p^\rmi(k)$ is the \emph{index} of $k$ in $p$.
 A permutation $\pi$ in $S_n$ that avoids the pattern $p$ must also avoid the
 bivincular pattern
 \[
 	\bivincs{8}{1/1/\cdot, 2/2/\cdot, 3/\cdot/\cdot, 4/\cdot/k,
	            5/\cdot/1, 6/\cdot/\cdot, 7/\cdot/\cdot, 8/k/\cdot}{}{2,3,4,5,6}
	=
 	(\dotsm k1 \dotsm, \tom, \{2,3,\dotsc,k-2\}).
 \]
\end{corollary}

By repeatedly applying the operations of inverse, reverse and complement
we can generate six other corollaries. We will not need them here.

\begin{example}
 Let's look at some simple applications:
 \begin{enumerate}
 
 \item Consider the bivincular pattern $p_1 = \vinc{4}{1/3, 2/1, 3/4, 4/2}{2}$. Proposition
 \ref{prop:bivinc-bivinc} shows that a permutation $\pi$ that avoids $p_1$ must
 also avoid $\bivinc{4}{1/3, 2/1, 3/4, 4/2}{}{2}$. In fact, the converse can be
 shown to be true, by taking inverses and applying the proposition.
 We will say more about the pattern $p_1$ in Example \ref{exampl:Baxter}.
 
 \item Consider the bivincular pattern $p_2 = \vinc{5}{1/3, 2/1, 3/5, 4/2, 5/4}{2}$.
 The proposition shows that a permutation $\pi$ that avoids $p_2$ must also
 avoid $\bivinc{5}{1/3, 2/1, 3/5, 4/2, 5/4}{}{2,3}$. We
 will say more about the pattern $p_2$ in subsection \ref{subsec:Gorenstein-Bivinc}.
 
 \end{enumerate}
\end{example}

\begin{example} \label{exampl:Baxter}
 The \emph{Baxter permutations} were originally defined and studied in
 relation to the ``commuting function conjecture'' of Dyer, see Baxter~\cite{MR0184217},
 and were enumerated by Chung, Graham, Hoggatt and Kleiman~\cite{MR491652}.
 Gire~\cite{G93} showed that these permutations
 can also be described as those avoiding
 the barred patterns $41\overline{3}52$ and $25\overline{1}34$.
 It was then pointed out by Ouchterlony~\cite{O05} that this is equivalent to
 avoiding the vincular patterns $\vinc{4}{1/3, 2/1, 3/4, 4/2}{2}$ and
 $\vinc{4}{1/2, 2/4, 3/1, 4/3}{2}$.

 Similarly to what we did above we can show that the Baxter permutations can
 also be characterized as those avoiding the bivincular patterns
 $\bivinc{4}{1/3, 2/1, 3/4, 4/2}{}{2}$ and $\bivinc{4}{1/2, 2/4, 3/1, 4/3}{}{2}$,
 and this is essentially a translation of the description in \cite{MR491652}
 into bivincular patterns.
\end{example}

Finally, here is an example that shows the converse of Proposition
\ref{prop:bivinc-bivinc} is not true.
\begin{example}
 The permutation $\pi = 423165$ avoids the pattern
 $\bivinc{5}{1/2, 2/3, 3/1, 4/5, 5/4}{}{2,3}$ but contains the pattern
 $\vinc{5}{1/2, 2/3, 3/1, 4/5, 5/4}{3}$, as the sub-word $23165$. 
\end{example}

\subsection{Gorenstein Schubert varieties in terms of bivincular patterns}
\label{subsec:Gorenstein-Bivinc}

Woo and Yong~\cite{MR2264071} classify those permutations $\pi$ that
correspond to Gorenstein Schubert varieties $X_\pi$. They do this using
embeddings of patterns with Bruhat restrictions, which we have described
above, and with a certain condition on the associated Grassmannian permutations
of $\pi$, which we will describe presently:

First, a \emph{descent} in a permutation $\pi$ is an integer $d$ such that
$\pi(d) > \pi(d+1)$, or equivalently, the index of the first letter in
an occurrence of the pattern $\vinc{2}{1/2, 2/1}{1}$.
A \emph{Grassmannian permutation} is a permutation
with a unique descent. Given any permutation $\pi$ we can associate
a Grassmannian permutation to each of its descents as follows: Given a
particular descent $d$ of $\pi$ we construct the sub-word $\gamma_d(\pi)$
by concatenating the right-to-left minima of the segment strictly to the
left of $d+1$ with the left-to-right maxima of the segment strictly to the
right of $d$. More intuitively we start with the descent $\pi(d)\pi(d+1)$
and enlarge it to the left by adding decreasing elements without creating another
descent and similarly enlarge it to the right by adding increasing elements
without creating another descent. We then denote the \emph{flattening}
(or \emph{standardization}) of $\gamma_d(\pi)$ by $\tilde{\gamma}_d(\pi)$,
which is the unique permutation whose letters are in the same relative order as
$\gamma_d(\pi)$.

\begin{example}
 Consider the permutation
 $\pi = 11 | 6 | 12 | 9 4 1 5 3 7 2 8 | 10$ where we have used the symbol
 $\sep$ to separate two-digit numbers from other numbers.
 For the descent at $d = 4$ we get $\gamma_4(\pi) = 6 9 4 5 7 8 | 10$
 and $\tilde{\gamma}_4(\pi) = 3 6 1 2 4 5 7$.
\end{example}

\noindent
Now, given a Grassmannian permutation $\rho$ in $S_n$ with its unique descent
at $d$ we construct its \emph{associated partition} $\lambda(\rho)$ as the
partition inside a bounding box $d \times (n-d)$, with $d$ rows and $n-d$ columns,
whose lower border is the lattice path that starts
at the lower left corner of the bounding box and whose $i$-th step,
for $i \in \dbrac{1,n}$, is vertical if $i$ is weakly to the left of
the position $d$, and horizontal otherwise. A corner of the lattice path is called
an \emph{inner corner} if it corresponds to a right turn on the path, otherwise it
is called an \emph{outer corner}. We are interested in the
\emph{inner corner distances} of this partition, that is, for every inner
corner we add its distance from the left side and the distance from the top
of the bounding box. If all these inner corner distances are the same then
the inner corners all lie on the same anti-diagonal.

In Theorem 1 of Woo and Yong~\cite{MR2264071} they show that a
permutation $\pi \in S_n$ corresponds to a Gorenstein Schubert variety $X_\pi$
if and only if
 
\begin{enumerate}
 
 \item \label{thm:WY06-cond1}
 for each descent $d$ of $\pi$, $\lambda(\tilde{\gamma}_d(\pi))$ has all of
 its inner corners on the same anti-diagonal; and
 
 \item \label{thm:WY06-cond2}
 the permutation $\pi$ avoids both $31524$ and
 $24153$ with Bruhat
 restrictions $\{ \brur{1}{5}, \brur{2}{3} \}$ and $\{ \brur{1}{5}, \brur{3}{4} \}$,
 respectively.
 
\end{enumerate}

Let's take a closer look at condition \ref{thm:WY06-cond2}: Proposition
\ref{prop:Bruhat-bivinc} below shows that
avoiding $31524$ with Bruhat restrictions
$\{ \brur{1}{5}, \brur{2}{3} \}$ is equivalent to avoiding the bivincular
pattern
\[
	\bivinc{5}{1/3,2/1,3/5,4/2,5/4}{2}{3} = (31524, \{2\}, \{3\} ).
\]
Similarly, avoiding
$24153$ with Bruhat restrictions $\{ \brur{1}{5}, \brur{3}{4} \}$
is equivalent to avoiding the bivincular
pattern
\[
	\bivinc{5}{1/2,2/4,3/1,4/5,5/3}{3}{2} = (24153, \{3\}, \{2\} ).
\]

\begin{proposition} \label{prop:Bruhat-bivinc}
 
 \begin{enumerate}
 
 \item \label{prop:Bruhat-bivinc-case1}
 Let $p$ be the pattern
 \[
 	\dotsm 1 k \dotsm
 \]
 in $S_k$. Let $j = p^\rmi(1)$ be the index of $1$ in $p$.
 A permutation $\pi$ in $S_n$ avoids $p$ with Bruhat restriction
 $\brur{j}{j+1}$ if and only if $\pi$ avoids the vincular pattern
 \[
 	\vinc{8}{1/\cdot, 2/\cdot, 3/\cdot, 4/1, 5/k, 6/\cdot, 7/\cdot, 8/\cdot}{4}
	=
	(\dotsm 1 k \dotsm, \{j\}, \tom).
 \]
 
 \item \label{prop:Bruhat-bivinc-case2}
 Let $\ell \in \dbrac{1,k-1}$ and $p$ be the pattern
 \[
 	\ell \dotsm (\ell+1)
 \]
 in $S_k$.
 A permutation $\pi$ in $S_n$ avoids $p$ with Bruhat restriction
 $\brur{1}{k}$ if and only if $\pi$ avoids the bivincular pattern
 \[
 	\bivincs{11}{1/1/\ell, 2/\cdot/\cdot, 3/\cdot/\cdot, 4/\ell/\cdot,
	 5/{\ell}/\cdot, 6/{+}/\cdot, 7/1/\cdot, 8/\cdot/\cdot, 9/\cdot/\ell, 10/\cdot/{+}, 11/k/1}
	 {}{4,5,6}
 	=
	(\ell \dotsm (\ell+1), \tom, \{ \ell \}).
 \]
 
 \end{enumerate}
 
\end{proposition}

\begin{proof}
 We consider each case separately.
 \begin{enumerate}
 
 \item Assume $\pi$ contains the vincular pattern mentioned. Then it clearly also
 contains an embedding satisfying the Bruhat restriction.
 
 Conversely assume $\pi$ contains an embedding satisfying the Bruhat restriction.
 If $1_\pi$ and $k_\pi$ are adjacent then we are done. Otherwise look at
 the element immediately to the right of $1_\pi$. This element must be either larger
 than $k_\pi$, in which case we can replace $k_\pi$ by it and are done, or smaller,
 in which case we replace $1_\pi$ by it, and repeat. 
 
 \item Assume $\pi$ contains the bivincular pattern mentioned. Then it clearly also
 contains an embedding satisfying the Bruhat restriction.
 
 Conversely assume $\pi$ contains an embedding satisfying the Bruhat restriction.
 If $(\ell+1)_\pi = \ell_\pi + 1$ then we are done. Otherwise consider the element
 $\ell_\pi + 1$. It must either be to the right of $(\ell+1)_\pi$ or to the left
 of $\ell_\pi$. In the first case we can replace $(\ell+1)_\pi$ by $\ell_\pi + 1$ and be
 done. In the second case replace $\ell_\pi$ with $\ell_\pi + 1$ and repeat. \qedhere
 \end{enumerate}
\end{proof}

As a consequence we get:

\begin{corollary}
 A permutation $\pi$ in $S_n$ avoids
 \[
 	\dotsm 1 k \dotsm,\ \brur{j}{j+1},
 \]
 where $j$ is the index of $1$, if and only if the inverse $\pi^\rmi$ avoids
 \[
 j \dotsm (j+1),\ \brur{1}{k}.
 \]
\end{corollary}

Note that we could have proved the statement of the corollary without
going through bivincular patterns and then used that to prove part
\ref{prop:Bruhat-bivinc-case2} of Proposition \ref{prop:Bruhat-bivinc},
as part \ref{prop:Bruhat-bivinc-case2} is the inverse statement of
the statement in part \ref{prop:Bruhat-bivinc-case1}.

Translating condition \ref{thm:WY06-cond1} of Theorem 1 of Woo and Yong~\cite{MR2264071}
into bivincular patterns is a bit more work. The failure of this condition is easily seen to
be equivalent to some partition $\lambda$ of an associated Grassmannian
permutation $\tilde{\gamma}_d(\pi)$
having an outer corner that is either ``too wide'' or ``too deep''.
More precisely, given a Grassmannian permutation $\rho$
and an outer corner of $\lambda(\rho)$,
we say that it is \emph{too wide} if the distance
upward from it to the next inner corner is smaller than the distance to the
left from it to the next inner corner. Conversely we say that an outer
corner is \emph{too deep} if the distance
upward from it to the next inner corner is larger than the distance to the
left from it to the next inner corner. We say that an outer corner
is \emph{unbalanced} if it is either too wide or too deep. We say that
an outer corner is \emph{balanced} if it is not \emph{unbalanced}.

If a permutation has an associated Grassmannian permutation with an outer corner
that is too wide we say
that the permutation itself is \emph{too wide} and similarly for
\emph{too deep}. If the permutation is either too wide or too deep we
say that it is \emph{unbalanced}, otherwise it is \emph{balanced}.
It is time to see some examples.

\begin{example} \label{ex:partitions}
See Figure \ref{fig:partitions} for drawings of the partitions below.
 \begin{enumerate}
 
 \item Consider the permutation $\rho = 14235$ with a unique descent at
 $d = 2$. It corresponds to the partition $(2) \subseteq 2 \times 3$
 and has just one outer corner. This outer corner is too wide.
 
 \item Consider the permutation $\rho = 13425$ with a unique descent at
 $d = 3$. It corresponds to the partition $(1,1) \subseteq 3 \times 2$
 and has just one outer corner. This outer corner is too deep.
 
 \item Consider the permutation $\rho = 1 3 4 8 9 2 5 6 7 | 10$ with a
 unique descent at $d = 5$. It corresponds to the partition
 $(4,4,1,1) \subseteq 5 \times 5$
 and has two outer corners. The first outer corner is too deep
 and the second is too wide.
 
 \item Consider the permutation $\rho = 1 3 6 7  2 4 5 8$ with a
 unique descent at $d = 4$. It corresponds to the partition
 $(3,3,1) \subseteq 4 \times 4$
 and has two outer corners that are both balanced.
 
 \end{enumerate}
 
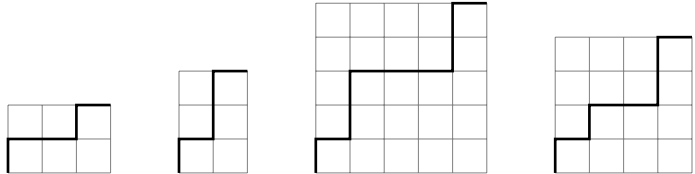
\begin{figure}[ht]
\begin{tikzpicture}[scale = 0.45, dot/.style = {fill,draw,circle, minimum size = 1pt}]

\draw[help lines] (0,0) grid (3,2);

\draw[line width = 1pt]
(0,0) -- (0,1) -- (1,1) -- (2,1) -- (2,2) -- (3,2) ;

\begin{scope}[xshift = 5cm]
\draw[help lines] (0,0) grid (2,3);

\draw[line width = 1pt]
(0,0) -- (0,1) -- (1,1) -- (1,2) -- (1,3) -- (2,3) ;
\end{scope}

\begin{scope}[xshift = 9cm]
\draw[help lines] (0,0) grid (5,5);

\draw[line width = 1pt]
(0,0) -- (0,1) -- (1,1) -- (1,2) -- (1,3) -- (2,3) -- (3,3) -- (4,3) -- (4,4) -- (4,5)-- (5,5) ;
\end{scope}

\begin{scope}[xshift = 16cm]
\draw[help lines] (0,0) grid (4,4);

\draw[line width = 1pt]
(0,0) -- (0,1) -- (1,1) -- (1,2) -- (2,2) -- (3,2) -- (3,3) -- (3,4) -- (4,4) ;
\end{scope}

\end{tikzpicture}

\caption{The associated partitions of the permutations in Example \ref{ex:partitions}}
\label{fig:partitions}
\end{figure}

\end{example}

We now show how these properties of Grassmannian permutations can be detected with
bivincular patterns.

\begin{lemma}
\label{lem:Grassmann-perm}
 Let $\rho$ be a Grassmannian permutation.
 
 \begin{enumerate}
 
 \item \label{lem:Grassmann-perm-too-wide}
 The permutation $\rho$ is too wide if and only if it contains at least one
 of the bivincular patterns from the infinite family
 \begin{align*}
 	\cF 
	&= \left\{
	\bivinc{5}{1/1, 2/4, 3/2, 4/3, 5/5}{}{2,3,4},
	\bivinc{7}{1/1, 2/5, 3/6, 4/2, 5/3, 6/4, 7/7}{}{2,3,4,5,6},
	\bivinc{9}{1/1, 2/6, 3/7, 4/8, 5/2, 6/3, 7/4, 8/5, 9/9}{}{2,3,4,5,6,7,8}, \dotsc
	\right\}.
 \end{align*}
 The general member of this family is of the form
 \[
 	\bivincs{11}{1/1/1, 2/2/\ell, 3/\cdot/{+}, 4/\cdot/1, 5/\cdot/\cdot,
	6/\cdot/\cdot, 7/\cdot/2, 8/\cdot/\cdot, 9/\cdot/\cdot, 10/\cdot/\ell, 11/k/k}
	{}{2,3,4,5,6,7,8,9,10},
 \]
 where $\ell = (k+1)/2$, and $k$ is odd.
 
 \item \label{lem:Grassmann-perm-too-deep}
 The permutation $\rho$ is too deep if and only if it contains at least one
 of the bivincular patterns from the infinite family
  \begin{align*}
 	\cG
	&= \left\{
	\bivinc{5}{1/1, 2/3, 3/4, 4/2, 5/5}{}{1,2,3},
	\bivinc{7}{1/1, 2/4, 3/5, 4/6, 5/2, 6/3, 7/7}{}{1,2,3,4,5},
	\bivinc{9}{1/1, 2/5, 3/6, 4/7, 5/8, 6/2, 7/3, 8/4, 9/9}{}{1,2,3,4,5,6,7}, \dotsc
	\right\}.
 \end{align*}
  The general member of this family is of the form
 \[
 	\bivincs{11}{1/1/1, 2/2/\ell, 3/\cdot/{+}, 4/\cdot/1, 5/\cdot/\cdot,
	6/\cdot/\cdot, 7/\cdot/2, 8/\cdot/\cdot, 9/\cdot/\cdot, 10/\cdot/\ell, 11/k/k}
	{}{1,2,3,4,5,6,7,8,9},
 \]
  where $\ell = (k-1)/2$, and $k$ is odd.
 
 \end{enumerate}
\end{lemma}

Note that these two infinite families can be obtained from
one another by reverse complement.

\begin{proof}
 We only consider part \ref{lem:Grassmann-perm-too-wide}, as part
 \ref{lem:Grassmann-perm-too-deep} is proved analogously. Assume
 that $\rho$ is a Grassmannian permutation that is too wide, so it
 has an outer corner that is too wide. Let $\ell$ be the distance
 from this outer corner to the next inner corner above. Then
 the distance from this outer corner to the next inner corner to
 the left is at least $\ell + 1$. This allows us to construct
 an increasing sequence $t$ of length $\ell$ in $\rho$,
 starting at a distance at least two to the right of the descent.
 We can also choose $t$ so that every element in it is
 adjacent both in location and values. Similarly we can construct
 an increasing sequence $s$ of length $\ell$ in $\rho$,
 located strictly to the left of the descent.
 We can also choose $s$ so that every element in it is
 adjacent both in location and values.
 This produces the required member of the family $\cF$.
 
 Conversely, assume $\rho$ contains the $i$-th member of the family
 $\cF$, the pattern
  \[
 	\bivincs{11}{1/1/1, 2/2/\ell, 3/\cdot/{+}, 4/\cdot/1, 5/\cdot/\cdot,
	6/\cdot/\cdot, 7/\cdot/2, 8/\cdot/\cdot, 9/\cdot/\cdot, 10/\cdot/\ell, 11/k/k}
	{}{2,3,4,5,6,7,8,9,10},
 \]
 where $k = 2i + 3$. Then the occurrence of the pattern
 corresponds to an outer corner in the partition of $\rho$ of width
 $\ell - 1$ and depth $\ell - 2$.
\end{proof}

We have now shown that:

\begin{proposition}
 A permutation $\pi$ is balanced if and only if every associated Grassmannian
 permutation avoids every bivincular pattern in the two infinite families
 $\cF$ and $\cG$ in Lemma \ref{lem:Grassmann-perm}.
\end{proposition}

This gives us:

\begin{theorem}
\label{thm:Ulfarsson-Gorenstein}
 Let $\pi \in S_n$. The Schubert variety $X_\pi$ is Gorenstein if and only if
 
 \begin{enumerate}
 
 \item $\pi$ is balanced; and
 
 \item
 the permutation $\pi$ avoids the bivincular patterns
 \[
 \bivinc{5}{1/3,2/1,3/5,4/2,5/4}{2}{3} \textrm{ and }
 \bivinc{5}{1/2,2/4,3/1,4/5,5/3}{3}{2}.
 \]
 
 \end{enumerate}
\end{theorem}

With these descriptions of factorial and Gorenstein varieties it is simple
to show that smoothness implies factoriality, which implies Gorensteinness:
If a variety is not factorial it contains either $\vinc{4}{1/2, 2/1, 3/4, 4/3}{2}$ or $1324$,
so it must contain either $2143$ and $1324$ and is therefore not smooth.
If a variety is not Gorenstein then at least one of the following are true,
\begin{enumerate}

\item $\pi$ has an associated Grassmannian permutation that contains
one of the bivincular patterns in the infinite families $\cF$ and
$\cG$ so it also contains $1324$ and is therefore not factorial.

\item $\pi$ contains $\bivinc{5}{1/3,2/1,3/5,4/2,5/4}{2}{3}$ or
$\bivinc{5}{1/2,2/4,3/1,4/5,5/3}{3}{2}$, so it also contains
$\vinc{4}{1/2, 2/1, 3/4, 4/3}{2}$ and is therefore not factorial.

\end{enumerate}

\section{Mesh patterns and marked mesh patterns} \label{sec:intmesh}

\subsection{Mesh patterns}
Br\"and\'en and Claesson~\cite{BC11} introduced a new type of pattern called a mesh pattern
and showed they generalize bivincular patterns and
(most) barred patterns. Here we recall their definition: A \emph{mesh pattern} is
a pair $(p,R)$ where $p$ is a permutation of rank $k$ and $R$ is a subset of
the square $\dbrac{0,k} \times \dbrac{0,k}$. An occurrence of that pattern in a permutation
$\pi$ is first of all an occurrence of $p$ in $\pi$ in the usual sense, that is,
a subset of the diagram $G(\pi) = \{ (i,\pi(i)) \sep 1 \leq i \leq n \}$. This occurrence
must also satisfy the restrictions determined by $R$, that is, there are order-preserving
injections $\alpha, \beta: \dbrac{1,k} \to \dbrac{1,n}$ such that if
$(i,j) \in R$ then $R_{ij} \cap G(\pi)$ is empty, where
\[
R_{ij} = \dbrac{\alpha(i)+1, \alpha(i+1)-1} \times
         \dbrac{\beta(j)+1,  \beta(j+1) -1};
\]
with $\alpha(0) = 0 = \beta(0)$ and $\alpha(k+1) = n+1 = \beta(k+1)$.
It is best to unwind this formal definition with a few examples.

\begin{example} \label{ex:markedmesh}
\begin{enumerate}
\item The mesh pattern
$(21, \{ (1,0),(1,1),(1,2) \})$ can be depicted as follows:
\[
\pattern{scale=1}{ 2 }{ 1/2, 2/1 }{1/0, 1/1, 1/2 }
\]
An occurrence of this mesh pattern in a permutation is an inversion
(an occurrence of the classical pattern $21$) with the additional requirement that there
is nothing in between the two elements in the occurrence. We usually refer
to this pattern as the vincular pattern $\vinc{2}{1/2,2/1}{1}$,
that is, a descent.

\item Now consider
the more complicated mesh pattern below.
\[
\pattern{scale=1}{ 2 }{ 1/1, 2/2 }{0/0, 2/0, 0/2, 2/2 }
\]
There are two occurrences of this mesh pattern in the permutation
$\pi = 315426$, shown below.
\[
\imopattern{scale=1}{ 6 }{ 1/3, 2/1, 3/5, 4/4, 5/2, 6/6 }
{}{ 1/3, 6/6 }{ 
0/0, 0/1, 0/2,
0/6,
6/6,
6/0, 6/1, 6/2 }
\qquad
\imopattern{scale=1}{ 6 }{ 1/3, 2/1, 3/5, 4/4, 5/2, 6/6 }
{}{ 2/1, 6/6 }{ 
0/0, 1/0,
0/6, 1/6,
6/6,
6/0 }
\]
Every other occurrence of the classical pattern $12$ in $\pi$, \textit{e.g.},
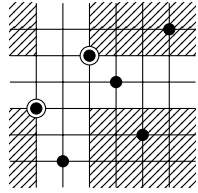
\[
\imopattern{scale=1}{ 6 }{ 1/3, 2/1, 3/5, 4/4, 5/2, 6/6 }
{}{ 1/3, 3/5 }{ 
0/0, 0/1, 0/2,
0/5, 0/6,
3/5, 4/5, 5/5, 6/5, 3/6, 4/6, 5/6, 6/6,
3/0, 4/0, 5/0, 6/0, 3/1, 4/1, 5/1, 6/1, 3/2, 4/2, 5/2, 6/2 }
\]
fails to satisfy the requirements given, since some of the shaded areas
will have a dot in them.

\item Dukes and Reifegerste~\cite{MR2628782} defined \emph{certified non-inversions} as occurrences
of $132$ that are neither part of $1432$ nor $1342$. Equivalently these are occurrences of the mesh
pattern below.
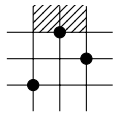
\[
\pattern{scale=1}{ 3 }{ 1/1, 2/3, 3/2 }{1/3,2/3}%
\]
\end{enumerate}
\end{example}

Br\"and\'en and Claesson~\cite{BC11} showed that a barred pattern with only one barred letter can
be written as a mesh pattern\footnote{This had been noticed earlier as well,
see unpublished work of Isaiah Lankham.}. The procedure is as follows: If $\pi(i)$
is the only barred letter in a barred pattern $\pi$ then the corresponding
mesh pattern is $(\pi', \{(i-1,\pi(i)-1)\})$ where $\pi'$ is the standardization
of $\pi$ after the removal of $\pi(i)$. For example
\[
1\overline{2}3 = \pattern{scale=1}{ 2 }{ 1/1, 2/2 }{1/1}.
\]
More general barred patterns can be also be translated into mesh patterns
as long as the barred letters are neither adjacent in locations nor in values.
The procedure is essentially the same as above, so for example the barred
pattern $63\overline{4}1\overline{2}5$ is contained in a permutation $\pi$
if and only if at least one of the mesh patterns
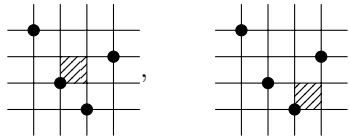
\[
\pattern{scale=1}{ 4 }{ 1/4, 2/2, 3/1, 4/3 }{2/2}, \qquad
\pattern{scale=1}{ 4 }{ 1/4, 2/2, 3/1, 4/3 }{3/1}
\]
is contained in $\pi$.

It is possible to classify simsun permutations by the avoidance of mesh patterns as follows:
A permutation in $S_n$ is \emph{simsun} if it contains no double descent ($\vinc{3}{1/3, 2/2, 3/1}{1,2}$)
in any of its restrictions to the interval $\dbrac{1,k}$ for some $k \leq n$.
For example the permutation $452613$ is not simsun since if we restrict it to $\dbrac{1,4}$
it becomes $4213$ which contains a double descent. It is now almost trivial to check that
a permutation is simsun if and only if it avoids the mesh pattern
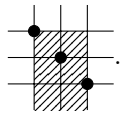
\[
\pattern{scale=1}{ 3 }{ 1/3, 2/2, 3/1 }{1/0, 1/1, 1/2, 2/0, 2/1, 2/2}.
\]
This was also noticed independently and simultaneously by Br\"and\'en and Claesson~\cite{BC11}. 

It is easy to see that bivincular patterns are special cases of mesh patterns. Adjacency
conditions on positions in the bivincular pattern become vertical strips, while adjacency
conditions on values become horizontal strips; $R$ is then the union of all the strips.
For example the bivincular
pattern $\bivinc{5}{1/3,2/1,3/5,4/2,5/4}{2}{3}$ from Theorem \ref{thm:Ulfarsson-Gorenstein}
corresponds to the mesh pattern
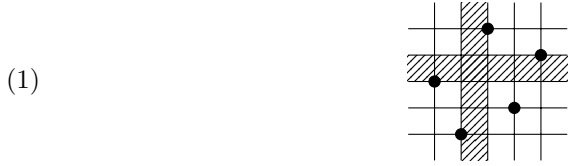
\begin{equation} \label{eq:transl}
\pattern{scale=1}{ 5 }{ 1/3, 2/1, 3/5, 4/2, 5/4 }{0/3, 1/3, 2/3, 3/3, 4/3, 5/3, 2/0, 2/1, 2/2, 2/3, 2/4, 2/5}.
\end{equation}

We have seen in Proposition \ref{prop:Bruhat-bivinc} that some Bruhat-restricted patterns
can be turned into bivincular patterns. We now show that any Bruhat-restricted pattern can
be turned into a mesh pattern.

Given a pattern $p$ with one Bruhat restriction $\brur{a}{b}$ first note
that this means that $a < b$ and $p(a) < p(b)$. Then a permutation contains $p$ with
the restriction if and only if it contains the mesh pattern
$(p, R)$, where $R$ consists of all the squares in the region with corners
$(a,\pi(a))$, $(b,\pi(a))$, $(b,\pi(b)$, $(a,\pi(b))$.
For example the pattern $31524$, $\brur{1}{5}$ corresponds to
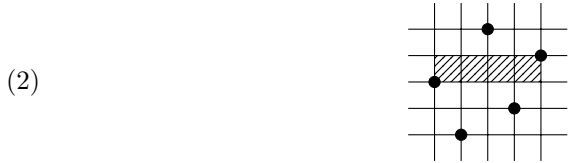
\begin{equation} \label{eq:transl2}
\pattern{scale=1}{ 5 }{ 1/3, 2/1, 3/5, 4/2, 5/4 }{1/3, 2/3, 3/3, 4/3}.
\end{equation}
Given a pattern $p$ with multiple Bruhat restrictions we superpose the
mesh patterns we get for each individual restriction. For example the
pattern $31524$, $\brur{1}{5}$, $\brur{3}{4}$, which is one of the patterns
that determines whether a permutation is Gorenstein or not, corresponds to
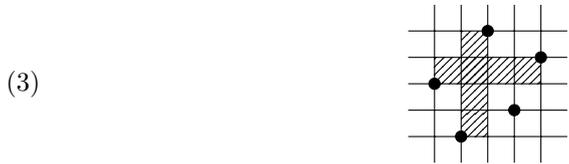
\begin{equation} \label{eq:transl3}
\pattern{scale=1}{ 5 }{ 1/3, 2/1, 3/5, 4/2, 5/4 }{1/3, 2/3, 3/3, 4/3, 2/1, 2/2, 2/3, 2/4}.
\end{equation}
Recall that we had already shown (Proposition \ref{prop:Bruhat-bivinc})
that this Bruhat-restricted pattern corresponds to the bivincular pattern
\eqref{eq:transl}. It is easy to see directly that the mesh pattern
\eqref{eq:transl3} is equivalent to the bivincular pattern \eqref{eq:transl},
in terms of being contained/avoided by a permutation.

It is now possible to translate Theorem
\ref{thm:Ulfarsson-Gorenstein} into mesh patterns, and completely get
rid of the middle step of considering the Grassmannian subpermutations.
But this was essentially done in Woo and Yong~\cite{MR2422304} using interval patterns, which
we now show to be special cases of mesh patterns.

\subsection{Interval patterns}
Woo and Yong~\cite{MR2422304} defined the avoidance of \emph{interval patterns}
as a generalization of Bruhat-restricted patterns.
We recall the definition here, with the modification that we reverse the usual
Bruhat order on $S_n$. We do this so the definition can be directly compared
with the definition of Bruhat-restricted avoidance. The (\emph{reversed})
\emph{Bruhat order} on $S_n$ is the partial order defined by
$\rho < \pi$ if $\pi$ can be obtained from $\rho$ by composing with a transposition
and $\pi$ has more non-inversions than $\rho$. Recall that a \emph{non-inversion}
is an occurrence of the classical pattern $12$; we let $\ell(\pi)$ denote the number of
non-inversions in $\pi$. Now we say that a permutation
$\pi$ \emph{contains the interval} $[p,q]$ if there exists a permutation $\rho \leq \pi$
and a common embedding of $p$ into $\rho$ and $q$ into $\pi$ such that the entries
outside of the embedding agree and the posets $[p,q]$, $[\rho,\pi]$ are isomorphic. 

\begin{example}
The interval pattern $[41523,31524]$
corresponds to the Bruhat-restricted pattern $31524$, $\brur{1}{5}$, shown
as a mesh pattern above, \eqref{eq:transl2}; and $[45123,31524]$ to $31524$, $\brur{1}{5}$,
$\brur{3}{4}$ also shown above, \eqref{eq:transl3}.
\end{example}

To show that any interval pattern can be turned into a mesh pattern we need a preliminary
definition: Given a permutation $\pi$ of rank $n$ and integers $j,k \in \dbrac{1,n+1}$,
we define a new permutation $\pi \oplus_j k$ of rank $n+1$ as follows:
\[
(\pi \oplus_j k) (\ell) =
\begin{cases}
\pi(\ell)          & \text{ if $\ell < j$ and $\pi(\ell) < k$}, \\
\pi(\ell)+1      & \text{ if $\ell < j$ and $\pi(\ell) \geq k$}, \\
k               & \text{ if $\ell = j$}, \\
\pi(\ell-1)          & \text{ if $\ell \geq j$ and $\pi(\ell) < k$}, \\
\pi(\ell-1)+1      & \text{ if $\ell \geq j$ and $\pi(\ell) \geq k$}.
\end{cases}
\]
For example $34125 \oplus_3 4 = 354126$.

\begin{lemma}
A permutation $\pi$ contains an interval pattern $[p,q]$ if and only
if it contains the mesh pattern $(q,R)$ where $R$ consists of boxes
$(i,j)$ such that
\[
\ell(q) - \ell(p) \neq \ell(q\oplus_i j) - \ell(p\oplus_i j)
\]
\end{lemma}
\begin{proof}
This lemma is a direct corollary of Lemma 2.1 in Woo and Yong~\cite{MR2422304}
which states that an embedding $\Phi$ of $[p,q]$ into $[\rho,\pi]$
is an interval pattern embedding if and only if $\ell(q)-\ell(p) = \ell(\pi) - \ell(\rho)$.
\end{proof}

It should be noted that although the general definition of a mesh pattern did not
exist many authors had drawn diagrams analogous to the diagrams we have been drawing
for mesh patterns, see e.g., \cite{MR1990570}, \cite{MR2376109}.

\begin{example}
To realize the interval $[53241,32154]$ as a mesh pattern we start
by drawing $53241$ with white dots and $32154$ with black dots into
the same diagram and consider the boxes $(i,j)$ that satisfy the condition
in the lemma above.
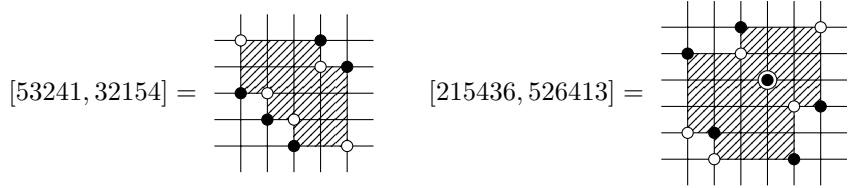
\[
[53241,32154] =
\impattern{scale=1}{ 5 }{ 1/3, 2/2, 3/1, 4/5, 5/4 }
{ 1/5, 2/3, 3/2, 4/4, 5/1 }{ 1/3, 1/4, 2/2, 2/3, 2/4, 3/1, 3/2, 3/3, 3/4, 4/1, 4/2, 4/3 }
\qquad
[215436, 526413] =
\imopattern{scale=1}{ 6 }{ 1/5, 2/2, 3/6, 4/4, 5/1, 6/3 }
{ 1/2, 2/1, 3/5, 5/3, 6/6 }{ 4/4 }{ 1/2, 1/3, 1/4, 2/1, 2/2, 2/3, 2/4, 3/1, 3/2, 3/3, 3/4, 3/5, 4/1, 4/2, 4/3, 4/4, 4/5, 5/3, 5/4, 5/5 }
\]
\end{example}

In Theorem 6.6 of Woo and Yong~\cite{MR2422304} they show that a
permutation $\pi \in S_n$ corresponds to a Gorenstein Schubert variety $X_\pi$
if and only if $\pi$ avoids intervals of the form
 
\begin{enumerate}
 
 \item \label{thm:WY08-cond1}
 $g_{a,b} = [(a+2) \dotsm (a+b+2) 1 \dotsm a (a+1), 1 (a+2) \dotsm (a+b+1) 2 \dotsm a (a+1) (a+b+2)]$
 for all integers $a,b > 0$ such that $a \neq b$,
 
 \item \label{thm:WY08-cond2}
 $h_{a,b} = [(a+4) \dotsm (a+b+4) (a+2) (a+3) 1 \dotsm (a+1), (a+2) (a+4) \dotsm (a+b+3) 1 (a+b+4) 2 \dotsm (a+1) (a+2)]$ 
 for all integers $a,b \geq 0$ such that either $a > 0$ or $b > 0$.
 
\end{enumerate}

\begin{figure}
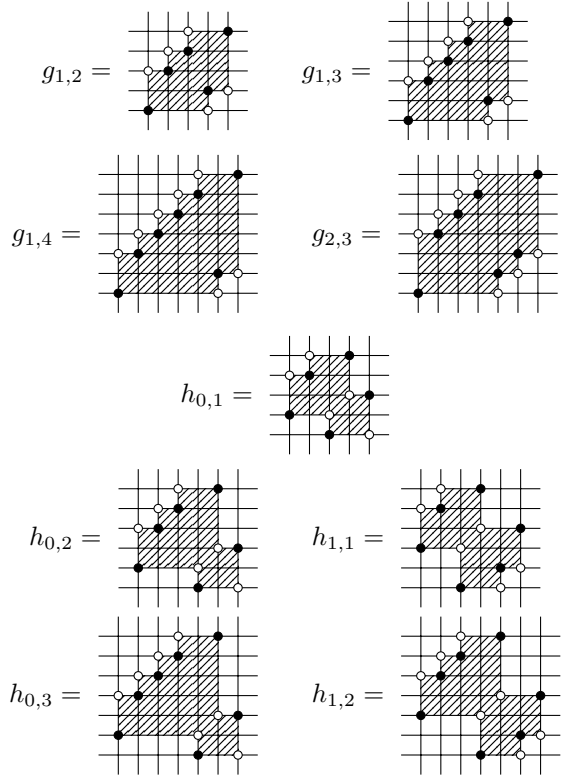

\[
g_{1,2} =
\impattern{scale=0.75}{ 5 }{ 1/1, 2/3, 3/4, 4/2, 5/5 }
{ 1/3, 2/4, 3/5, 4/1, 5/2 }{1/1, 1/2, 2/1, 2/2, 2/3, 3/1, 3/2, 3/3, 3/4, 4/2, 4/3, 4/4}
\qquad
g_{1,3} =
\impattern{scale=0.75}{ 6 }{ 1/1, 2/3, 3/4, 4/5, 5/2, 6/6 }
{ 1/3, 2/4, 3/5, 4/6, 5/1, 6/2 }{1/1, 1/2, 2/1, 2/2, 2/3, 3/1, 3/2, 3/3, 3/4, 4/1, 4/2, 4/3, 4/4, 4/5, 5/2, 5/3, 5/4, 5/5}
\]
\[
g_{1,4} =
\impattern{scale=0.75}{ 7 }{ 1/1, 2/3, 3/4, 4/5, 5/6, 6/2, 7/7 }
{ 1/3, 2/4, 3/5, 4/6, 5/7, 6/1, 7/2 }{1/1, 1/2, 2/1, 2/2, 2/3, 3/1, 3/2, 3/3, 3/4, 4/1, 4/2, 4/3, 4/4, 4/5, 5/1, 5/2, 5/3, 5/4, 5/5, 5/6, 6/2, 6/3, 6/4, 6/5, 6/6}
\qquad
g_{2,3} =
\impattern{scale=0.75}{ 7 }{ 1/1, 2/4, 3/5, 4/6, 5/2, 6/3, 7/7 }
{ 1/4, 2/5, 3/6, 4/7, 5/1, 6/2, 7/3 }{1/1, 1/2, 1/3, 2/1, 2/2, 2/3, 2/4, 3/1, 3/2, 3/3, 3/4, 3/5, 4/1, 4/2, 4/3, 4/4, 4/5, 4/6, 5/2, 5/3, 5/4, 5/5, 5/6, 6/3, 6/4, 6/5, 6/6}
\]

\[
h_{0,1} =
\impattern{scale=0.75}{ 5 }{ 1/2, 2/4, 3/1, 4/5, 5/3 }
{ 1/4, 2/5, 3/2, 4/3, 5/1 }{1/2, 1/3, 2/2, 2/3, 2/4, 3/1, 3/2, 3/3, 3/4, 4/1, 4/2}
\]
\[
h_{0,2} =
\impattern{scale=0.75}{ 6 }{ 1/2, 2/4, 3/5, 4/1, 5/6, 6/3 }
{ 1/4, 2/5, 3/6, 4/2, 5/3, 6/1 }{1/2, 1/3, 2/2, 2/3, 2/4, 3/2, 3/3, 3/4, 3/5, 4/1, 4/2, 4/3, 4/4, 4/5, 5/1, 5/2}
\qquad
h_{1,1} =
\impattern{scale=0.75}{ 6 }{ 1/3, 2/5, 3/1, 4/6, 5/2, 6/4 }
{ 1/5, 2/6, 3/3, 4/4, 5/1, 6/2 }{1/3, 1/4, 2/3, 2/4, 2/5, 3/1, 3/2, 3/3, 3/4, 3/5, 4/1, 4/2, 4/3, 5/2, 5/3}
\]
\[
h_{0,3} =
\impattern{scale=0.75}{ 7 }{ 1/2, 2/4, 3/5, 4/6, 5/1, 6/7, 7/3 }
{ 1/4, 2/5, 3/6, 4/7, 5/2, 6/3, 7/1 }{1/2, 1/3, 2/2, 2/3, 2/4, 3/2, 3/3, 3/4, 3/5, 4/2, 4/3, 4/4, 4/5, 4/6, 5/1, 5/2, 5/3, 5/4, 5/5, 5/6, 6/1, 6/2}
\qquad
h_{1,2} =
\impattern{scale=0.75}{ 7 }{ 1/3, 2/5, 3/6, 4/1, 5/7, 6/2, 7/4 }
{ 1/5, 2/6, 3/7, 4/3, 5/4, 6/1, 7/2 }{1/3, 1/4, 2/3, 2/4, 2/5, 3/3, 3/4, 3/5, 3/6, 4/1, 4/2, 4/3, 4/4, 4/5, 4/6, 5/1, 5/2, 5/3, 6/2, 6/3}
\]
 \caption{The first few patterns in Theorem 6.6 of Woo and Yong~\cite{MR2422304} shown as mesh patterns (patterns
that are the reverse complement of a pattern that has already appeared are omitted)}
\label{fig:WooYongGor}
\end{figure}

\noindent
See Figure \ref{fig:WooYongGor} for some patterns appearing in these two lists.

\subsection{Marked mesh patterns, DBI- and Gorenstein varieties}%
We introduce a generalization of mesh patterns which we call
\emph{marked mesh patterns} and use them to give an alternative description of Schubert
varieties defined by inclusions, Gorenstein Schubert varieties, $123$-hexagon
avoiding permutations, Dumont permutations and cycles in permutations.

\begin{definition}
A \emph{marked mesh pattern} $(p,\mathcal{C})$ of rank $k$ consists of a classical pattern
$p$ of rank $k$ and a collection $\mathcal{C}$ which contains pairs $(C, \square j)$
where $C$ is a subset of the square $\dbrac{0,k} \times \dbrac{0,k}$, $j$ is a non-negative
integer and $\square$ is one of the symbols $\leq, =, \geq$. An occurrence of $(p,\mathcal{C})$ in a permutation $\pi$
is first of all an occurrence of $p$ in $\pi$ in the usual sense, that is, a subset
of the diagram $G(\pi) = \{ (i,\pi(i)) \sep 1 \leq i \leq n \}$. This occurrence
must also satisfy the restrictions determined by $\mathcal{C}$, that is, there are order-preserving
injections $\alpha, \beta: \dbrac{1,k} \to \dbrac{1,n}$ such that for each pair $(C, \square j)$
we have
\[
\#(C' \cap G(\pi))\ \square\ j,
\]
where
\[
C' = \bigcup_{(i,j) \in C} R_{ij}.
\]
As above, 
$R_{ij} = \dbrac{\alpha(i)+1, \alpha(i+1)-1} \times \dbrac{\beta(j)+1,  \beta(j+1) -1}$,
with $\alpha(0) = 0 = \beta(0)$ and $\alpha(k+1) = n+1 = \beta(k+1)$.
\end{definition}

\noindent
Since regions of the form $(C,\geq j)$ are the most common we also write them more simply as $(C,j)$. 

\begin{example}
\begin{enumerate}

\item Every mesh pattern $(p,R)$ is an example of a marked mesh pattern, we just define
$\mathcal{C} = \{(R,{=}0)\}$. For example, here is the mesh pattern that identifies
the simsun permutations, written as a marked mesh pattern.
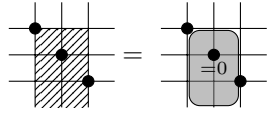
\[
\pattern{scale=1}{ 3 }{ 1/3, 2/2, 3/1 }{1/0, 1/1, 1/2, 2/0, 2/1, 2/2}
=
\patternsbmm{scale=1}{ 3 }{ 1/3, 2/2, 3/1 }{}{1/0/3/3/k}{1/0/3/3/{=0}}
\]

\item Consider the marked mesh pattern
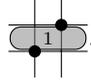
\[
\patternsbm{scale=1}{ 2 }{ 1/1, 2/2 }{}{0/1/3/2/1}.
\]
If a permutation $\pi$ contains it, it has an occurrence of the classical pattern $12$ where
there is at least one element $x$ in the permutation with the property that $1_\pi < x < 2_\pi$.
This is equivalent to saying that $\pi$ contains at least one of the classical patterns $213, 123, 132$.

\item A fixed point of a permutation is an occurrence of the marked mesh pattern
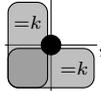
\[
\patternsbmm{scale=1.75}{ 1 }{ 1/1 }{}{0/0/1/2/k, 0/0/2/1/k}{0/1/1/2/{=k}, 1/0/2/1/{=k}},
\]
for some integer $k \geq 0$.
This generalizes to occurrences of cycles in a permutation. For example, a $2$-cycle
is an occurrence the marked mesh pattern $(21,\mathcal{C})$ where
$\mathcal{C}$ consists of the four marked regions below
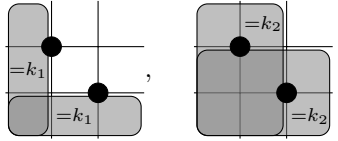
\[
\patternsbmm{scale=1.75}{ 2 }{1/2, 2/1}{}{0/0/1/3/k, 0/0/3/1/k}{0/1/1/2/{=k_1}, 1/0/2/1/{=k_1}}, \quad
\patternsbmm{scale=1.75}{ 2 }{1/2, 2/1}{}{0/0/2/3/k, 0/0/3/2/k}{1/2/2/3/{=k_2}, 2/0.5/3/0.5/{=k_2}},
\]
for some integers $k_2 > k_1 \geq 0$. These types of patterns can also be extended to include
unions of cycles and thus subsume the patterns defined in \mbox{McGovern}~\cite{MR2833467}.

\item \emph{Dumont permutations of the first kind}~\cite{0297.05004} are permutations of even rank with the property that every
even integer is followed by a smaller integer and every odd integer is either the last entry in the permutation
or is followed by a larger integer.
Therefore a permutation is a Dumont permutation of the first kind if and only if it avoids the marked mesh patterns
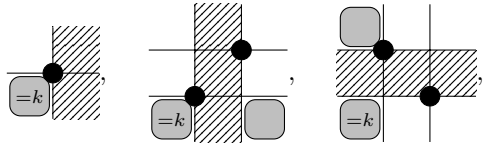
\[
\patternsbmm{scale=1.75}{ 1 }{1/1}{1/0,1/1}{0/0/1/1/k}{0/0/1/1/{=k}}, \quad
\patternsbmm{scale=1.75}{ 2 }{1/1,2/2}{1/0,1/1,1/2}{0/0/1/1/k, 2/0/3/1/k}{0/0/1/1/{=k}}, \quad
\patternsbmm{scale=1.75}{ 2 }{1/2,2/1}{0/1,1/1,2/1}{0/0/1/1/k, 0/2/1/3/k}{0/0/1/1/{=k}}, \quad
\]
where $k$ is an odd integer. Note that in the second pattern there is a single marked region $(\{(0,0),(2,0)\},{=}k)$,
consisting of two separated boxes. Similarly for the third pattern.
Dumont permutations of the second kind are also defined in~\cite{0297.05004}.
They can also be defined by the avoidance of the marked mesh patterns
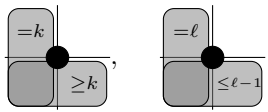
\[
\patternsbmm{scale=2}{ 1 }{ 1/1 }{}{0/0/1/2/k, 0/0/2/1/k}{0/1/1/2/{=k}, 1/0/2/1/{\geq k}}, \quad
\patternsbmm{scale=2}{ 1 }{ 1/1 }{}{0/0/1/2/k, 0/0/2/1/k}{0/1/1/2/{=\ell}, 1/0/2/1/{\scriptscriptstyle{\leq \ell-1}}}
\]
where $k$ is an odd integer and $\ell$ is an even integer.

\item Green and Losonczy~\cite{MR1980344} defined \emph{freely braided permutations} as those permutations
avoiding the classical patterns $3421$, $4231$, $4312$, $4321$. Equivalently, these are the permutations
avoiding the marked mesh pattern
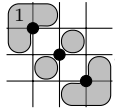
\[
\patternsbmfreelybraided{scale=1}{ 3 }{1/3,2/2,3/1}{}{1/1/2/2, 2/2/3/3}{0/3/1/4/1},
\]
marked with a single region consisting of $(2,0)$, $(3,0)$, $(1,1)$, $(3,1)$, $(0,2)$, $(2,2)$, $(3,0)$ and $(3,1)$.

\item Labelle, Leroux, Pergola and Pinzani~\cite{MR1887485} defined an \emph{inversion of the $j$-th kind}
in a permutation $\pi$ to be a pair of elements $\pi(s) > \pi(t)$, with $s < t$ and such that there
do not exist $j$ distinct indices $t+1 \leq t_1,t_2,\dotsc,t_j \leq n$ such that $\pi(t) > \pi(t_i)$
for $i = 1, \dotsc, j$.
Alternatively, an inversion of the $j$-th kind is an occurrence of the marked mesh pattern
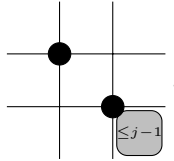
\[
\patternsbm{scale=2}{ 2 }{ 1/2, 2/1 }{}{2/0/3/1/{\scriptscriptstyle{\leq j-1}}}.
\]

\item Kitaev, see e.g.~\cite{MR2319869}, introduced \emph{partially ordered patterns} (\emph{POP})
as a generalization of vincular patterns. Some POPs can be written as marked mesh patterns, e.g.,
an occurrence of the POP $121$ in a permutation means the occurrence of either $231$ or $132$. Therefore
\[
121 = \patternsbm{scale=1}{ 1 }{ 1/1}{}{0/0/1/1/1, 1/0/2/1/1}.
\]

\item Hou and Mansour~\cite{1161.05301} studied permutations avoiding
$\vinc{4}{1/1, 2/{\ \square\ }, 3/{\, 2}, 4/3}{3}$ which are permutations
avoiding the marked mesh pattern
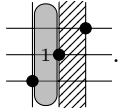
\[
\patternsbm{scale=1}{ 3 }{ 1/1, 2/2, 3/3 }{2/0, 2/1, 2/2, 2/3}{1/0/2/4/1}.
\]
\end{enumerate}
\end{example}

Gasharov and Reiner~\cite{MR1934291} defined Schubert varieties \emph{defined by inclusions} (or
just \emph{DBI}) and
characterized them with pattern avoidance of the patterns $24153$, $31524$, $426153$ and $1324$.
We now show how the first three of these patterns can be represented as two marked mesh patterns.

\begin{theorem} \label{thm:dbi}
Let $\pi \in S_n$. The Schubert variety $X_\pi$ is defined by inclusions if and only if
the permutation $\pi$ avoids the patterns
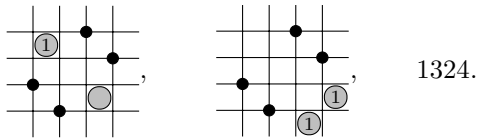
\[
\patternsbm{scale=1}{ 4 }{ 1/2, 2/1, 3/4, 4/3 }{}{3/1/4/2/{},1/3/2/4/1},
\qquad
\patternsbm{scale=1}{ 4 }{ 1/2, 2/1, 3/4, 4/3 }{}{4/1/5/2/1,3/0/4/1/1},
\qquad
1324.
\]
Where it should be noted that the first marked mesh pattern is marked with a single region,
$\{ (1,3), (3,1) \}$,
consisting of two boxes, and the number of dots in this region is at least $1$.
\end{theorem}
\begin{proof}
For the first marked mesh pattern note that $2143 \oplus_2 4 = 24153$ and $ 2143 \oplus_4 2 = 31524$.
For the second marked mesh pattern note that $(2143 \oplus_4 1) \oplus_6 4 = 426153$.
\end{proof}

We can also use these patterns to describe Gorenstein Schubert varieties:

\begin{theorem} \label{thm:Ulfarsson-Gorenstein2}
Let $\pi \in S_n$. The Schubert variety $X_\pi$ is Gorenstein if and only if it is balanced
and avoids the pattern
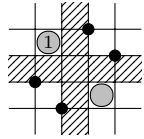
\[
\patternsbm{scale=1}{ 4 }{ 1/2, 2/1, 3/4, 4/3 }{2/0, 2/1, 2/2, 2/3, 2/4, 0/2, 1/2, 3/2, 4/2}{3/1/4/2/{},1/3/2/4/1}.
\]
\end{theorem}
\begin{proof}
Similar to the proof of Theorem \ref{thm:dbi}.
\end{proof}

In \cite{MR1826948} Billey and Warrington introduced \emph{$123$-hexagon avoiding}
permutations as permutations avoiding the classical patterns $123$,
$53281764$, $53218764$, $43281765$, $43218765$.\footnote{Actually they introduced \emph{$321$-hexagon avoiding}
permutations as permutations avoiding the reverse of the patterns listed here. We consider the reversed definition
to be compatible with what has appeared above.}
We now show how these four patterns can then be combined into one marked mesh pattern.

\begin{proposition} \label{prop:123hex}
A permutation $\pi$ is $123$-hexagon avoiding if and only if it avoids $123$ and
the marked mesh pattern
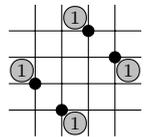
\[
\patternsbm{scale=1}{ 4 }{ 1/2, 2/1, 3/4, 4/3 }{}{2/4/3/5/{1}, 0/2/1/3/{1}, 4/2/5/3/{1}, 2/0/3/1/{1}}.
\]
\end{proposition}
\begin{proof}
The \lqq if\rqq\ part is easily verified. Now assume $\pi$ contains the marked mesh pattern.
Let $x$, $y$, $z$, $w$ correspond, respectively, to elements in the marked regions, read clockwise and
starting at the top. Let us assume that $\pi^\rmi(x) < \pi^\rmi(z)$ and $y < w$, as the other
cases are similar. Then $\pi$ contains the pattern $53281764$.
\end{proof}

Billey and Warrington also showed that $123$-hexagon avoiding permutations can be characterized
by the avoidance of $123$ and the avoidance of a hexagon in the \emph{heap} of the permutation.
See \cite{MR1826948}.

Tenner \cite{MR2333139} studied\footnote{Actually, permutations
that avoid the reverse of these patterns where considered.} permutations avoiding $123$ and $2143$ and showed
that a permutation $\pi$ avoids these two patterns if and only if it is \emph{boolean}
in the sense that the principal order ideal in strong Bruhat order $B(\pi)$ is \emph{Boolean} (that is,
isomorphic to the Boolean poset $B_r$ of subsets of $\dbrac{r}$ for some $r$).
So we immediately get that a Boolean permutation is $123$-hexagon avoiding.

The author is working with a coauthor on determining which patterns control
local complete intersections. The two marked mesh patterns that appear for Schubert
varieties defined by inclusions appear in the description along with one other marked mesh pattern.

\begin{figure}[ht]
 \begin{tikzpicture}
 [scale = 0.75, place/.style = {circle,draw = green!50,fill = green!20,thick,minimum size = 5pt},auto]
 
 \node[] at (0,6) (mm) {marked mesh};
 \node[] at (0,5) (m) {mesh};
 \node[] at (0,3) (biv) {bivincular};
 \node[] at (0,1) (vin) {vincular};
 \node[] at (0,0) (cl) {classical};
 
 \node[] at (-3,4) (1bar) {$1$-barred};
 \node[] at (3,4) (inter) {interval};
 \node[] at (3,2) (Bru) {Bruhat-restricted};
 
 \node[] at (-3,5) (McG) {McGovern};

 \draw[-] (mm) to (m);
 \draw[-] (mm) to (McG);
 \draw[-] (m) to (1bar);
 \draw[-] (m) to (biv);
 \draw[-] (m) to (inter);
 \draw[-] (inter) to (Bru);
 \draw[-] (biv) to (vin);
 \draw[-] (vin) to (cl);
 
 \end{tikzpicture}
 \caption{A diagram showing a hierarchy of permutation pattern definitions. $1$-barred refers
 to barred patterns with one bar}
 \label{fig:diagram}
\end{figure}
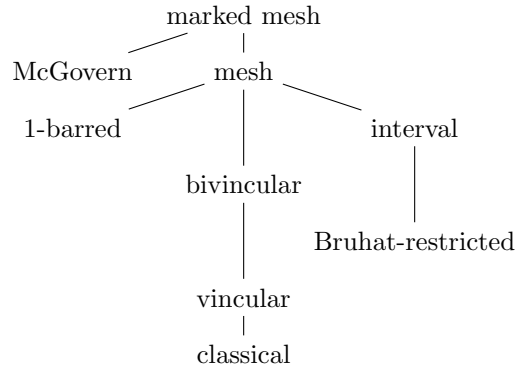

We end with a diagram in Figure \ref{fig:diagram} that shows which pattern definitions subsume which.

There are still other definitions of permutation patterns such as \emph{grid patterns}, defined by
Huczynska and Vatter~\cite{MR2240760}.
I do not know where they fit into the hierarchy in Figure \ref{fig:diagram}.

\section*{Acknowledgements}
The author thanks Einar Steingr\'imsson and Alexander Woo for many helpful conversations, and Sara Billey
for suggesting the title. The author also thanks an anonymous referee for many helpful comments.
This work is supported by grant no.\ 090038011 from the Icelandic Research Fund.

\bibliographystyle{amsplain}
\bibliography{RefsSchubert-Project}   

%
%

\end{document}